\documentclass[leqno]{article}

\usepackage{amsmath}
\usepackage{amssymb}
\usepackage{theorem}
\usepackage{mathtools}
\usepackage[utf8]{inputenc}
\usepackage[nottoc,numbib]{tocbibind}
\usepackage{multirow}
\usepackage{stmaryrd}
\usepackage{dutchcal}
\usepackage{textgreek}
\usepackage{url}

\usepackage[small,nohug,heads=vee]{diagrams}
\diagramstyle[lablestyle=\scriptstyle]

\theoremstyle{change}
\theorembodyfont{\itshape}
\theoremheaderfont{\normalfont\bfseries}
\newtheorem{thm}[equation]{Theorem}
\newtheorem{cor}[equation]{Corollary}
\newtheorem{prop}[equation]{Proposition}
\newtheorem{lem}[equation]{Lemma}

\theoremstyle{change}
\theorembodyfont{\normalfont}
\theoremheaderfont{\slshape}

\newtheorem{rem}[equation]{Remark}
\newtheorem{eg}[equation]{Example}

\newcommand{\Lg}{\mbox{$\mathfrak g$}}

\newcommand{\Lh}{\mbox{$\mathfrak h$}}
\newcommand{\Lk}{\mbox{$\mathfrak k$}}
\newcommand{\Lp}{\mbox{$\mathfrak p$}}
\newcommand{\La}{\mbox{$\mathfrak a$}}

\newcommand{\Lm}{\mbox{$\mathfrak m$}}

\newcommand{\Lq}{\mbox{$\mathfrak q$}}

\newcommand{\Ls}{\mbox{$\mathfrak s$}}

\newcommand{\Pf}{{\em Proof}. }
\newcommand{\EPf}{\hfill$\square$}

\newcommand{\ft}{\footnote}

\newcommand{\ad}{\mbox{$\mathrm{ad}$}}
\newcommand{\Ad}{\mbox{$\mathrm{Ad}$}}

\newcommand{\Z}{\mbox{$\mathbb Z$}}
\newcommand{\R}{\mbox{$\mathbb R$}}
\newcommand{\C}{\mbox{$\mathbb C$}}

\newcommand{\Ca}{\mbox{$\mathbb{O}$}}

\newcommand{\Hd}{\mbox{$\mathcal H$}}

\numberwithin{equation}{subsection}

\title{Topics in Polar Actions}
\author{Claudio Gorodski\thanks{\textsc{Instituto de Matem\'atica e Estat\'\i stica, Universidade de S\~ao Paulo, Rua do Mat\~ao, 1010, S\~ao Paulo, SP 05508-090, Brazil}, \textit{E-mail address}: \texttt{claudio.gorodski@usp.br}}}
\date{July 2022}

\begin{document}

\maketitle

\begin{abstract}
These are the notes for a series of lectures at the Institute of 
Geometry and Topology of the University of 
Stuttgart, Germany, in July 13-15,
2022. We wish to thank Uwe Semmelmann and Andreas Kollross
for the invitation to give these lectures.
We assume basic knowledge of isometric actions on Riemannian manifolds,
including the normal slice theorem and
the principal orbit type theorem. 
Lecture~1 introduces polar actions and culminates 
with Heintze, Liu and Olmos's argument to characterize 
them in terms of integrability of the distribution of normal 
spaces to the principal orbits. 
The other two lectures are devoted to two of Lytchak and Thorbergsson's 
results. In Lecture~2 we briefly review Riemannian orbifolds from the 
metric point of view, and explain their characterization of orbifold
points in the orbit space of a proper and isometric action in terms 
of polarity of the slice representation above. In 
Lecture~3 we present their proof of the fact that variationally complete 
actions in the sense of Bott and Samelson on non-negatively curved manifolds
are hyperpolar. The appendix contains explanations of some results used in the
lectures, namely: a criterion for
the polarity of isometric actions on symmetric spaces, 
a discussion of Cartan's and Hermann's criterions for the
existence of totally geodesic submanifolds, and a more or less self-contained 
derivation of Wilking's transversal Jacobi equation.

\end{abstract}

\section*{Introduction}

J. Dadok~\cite{dadok} considered orthogonal representations of compact Lie groups
with the property that there is a subspace meeting all orbits, and always
orthogonally. He pointed out that they resemble generalized
polar coordinates and introduced the name \emph{polar representations} for
them. In the same paper he proved a number of basic properties, and
especially classified polar representations up to orbit-equivalence\ft{Two isometric actions are said to be \emph{orbit-equivalent} 
if they have the same orbits, up to an isometric identification between the 
target spaces.} using highest weight theory.
As an aftermath, he realized that every polar representation of a
compact connected Lie group is orbit-equivalent to the isotropy
representation of a (Riemannian) symmetric 
space.

Polar actions had already been considered before, in one way or another.
One form of \'E. Cartan's maximal torus theorem states
that a maximal torus with a bi-invariant metric
meets every adjoint orbit in the group, orthogonally,
and similarly, the Lie algebra of the maximal torus meets every adjoint
orbit in the Lie algebra, orthogonally. A generalization of this
result applies to the (resp. linear) isotropy action of a symmetric space
and the (resp. Cartan subspaces) maximal flats.
L. Conlon~\cite{conlon}, building on work of Bott and Samelson
and Hermann, for an action of compact connected Lie group $K$ on a
complete Riemannian manifold
$M$, considered his so-called ``$K$-transversal domain'', namely,
a flat closed connected totaly geodesic embedded submanifold meeting
every $K$-orbit, and orthogonal at every point of intersection.
J. Szenthe, coming from a background on transformation groups,
and initially unaware of Conlon's work, studied the generalized
Weyl group for isometric actions of compact Lie groups admitting
``orthogonally transversal manifolds'', and proved that 
such submanifolds are automatically
totally geodesic. Palais and Terng~\cite{palais-terng}, initially unaware
of Dadok's and Conlon's works, defined a \emph{section} of an
isometric action of a compact Lie group on a connected complete
Riemannian manifold $M$ to be a
connected closed regularly embedded smooth submanifold $\Sigma$
of $M$ that meets all orbits orthogonally. They also note that
the compactness assumption on $G$ can be replaced by the hypothesis
that the action is proper, without substantial changes in the results.
They especially took the differential geometric viewpoint and emphasized
the relation to Riemannian geometry of submanifolds (more especifically,
isoparametric submanifolds, another area with important contributions
by Cartan). In the same paper, they mention applications
to invariant theory and calculus of variations. These
are further exposed in their book~\cite{palais-terng-book}.
Much more recently, later developments in the area of polar actions
are collected in~\cite{berndt-console-olmos}, which has also an extensive
list of bibliographic references. 

\emph{Singular Riemannian foliations}
form a class of foliations that generalize the foliations by orbits 
of an isometric group action. Much of the theory of polar actions 
has been generalized to ``singular Riemannian sections with sections'', 
or ``polar foliations'' as they are now called, but the lack
of group action causes some difficulties. We will not discuss them here  
and instead refer to~\cite{alexandrino-bettiol,radeschi,thorbergsson}
for discussions and references. 

\tableofcontents

\section{Lecture 1: Polar actions}

\subsection{Sections}

A proper isometric action of a Lie group $G$
on a complete Riemannian manifold $M$ is called
\emph{polar}\index{action!polar} if there
exists a complete connected immersed submanifold
$\Sigma$ of $M$ which intersects all the orbits 
and such that $\Sigma$ is perpendicular to every orbit it 
meets\ft{It is possible 
to consider a more general situation in which $\Sigma$ is replaced
by an isometric immersion~$\iota:\Sigma\to M$, non-necessarily injective.
In~\cite{berndt-console-olmos} such actions are called
\emph{locally polar}, and the term ``polar'' is reserved to those
actions with an embedded section. 
It is clear from the discussion below that $\iota$ can fail to be 
injective at~$p$ only if $\iota(p)$ is a singular point of the 
action. The brothers Alekseevsky~\cite{alekseevsky-alekseevsky} 
proved that an isometrically immersed section of dimension~$1$ 
must be injective, but the case of higher dimension remains open.
For simplicity of exposition, herein we restrict to injectively immersed
sections and just make some comments in the general case.}. 
Such a submanifold is called a \emph{section}.

A number of basic properties of polar actions 
is listed in Proposition~\ref{basic-properties-polar-action} below. First, we 
prove a related result about general proper isometric actions. 

\begin{lem}\label{exponential-normal-space}
Let $(G,M)$ be a proper isometric action.
Then, for every $p\in M$, the subset
$\exp_p(\nu_p(Gp))$ meets all the orbits of $G$.
\end{lem}

\Pf Fix an arbitrary orbit $N$ of $G$ and a point $q\in N$. Since the
action is proper, $N$ is a properly embedded, thus closed
submanifold of $M$. By completeness of~$M$, 
there exists a minimizing geodesic $\gamma:[0,1]\to M$ 
joining $q$ to $Gp$, so $\gamma(0)=q$ and $\gamma(1)\in Gp$.
Due the first variation of length formula, 
$\gamma$ is perpendicular to $Gp$ at $\gamma(1)$. 
Write $\gamma(1)=gp$ for some $g\in G$. Then 
$\tilde\gamma=g^{-1}\circ\gamma$ is a geodesic joining $g^{-1}q\in N$ 
to $p$, and it is perpendicular to $Gp$ at $p$. 
Therefore $g^{-1}q=\exp_pv$
where $v=-\tilde\gamma'(1)\in\nu_p(Gp)$.  This proves 
that $\exp_p(\nu_p(Gp))$ meets~$N$ at $g^{-1}q$. \EPf

\begin{prop}\label{basic-properties-polar-action}
Let $(G,M)$ be a polar action. Then:
\begin{enumerate}
\item If $\Sigma$ is a section of $(G,M)$ and $g\in G$, then $g\Sigma$ 
is a section of $(G,M)$. In other words, any $G$-translate 
of a section is a section.  
\item There exists a section of $(G,M)$ through every point of $M$.
\item The dimension of a section of $(G,M)$ equals the 
co-homogeneity of the action.  
\item Any section of $(G,M)$ contains an open and dense subset 
consisting of regular points of the action.
\item A section of $(G,M)$ is totally geodesic in $M$.
\item There exists a unique section of $(G,M)$ through 
a regular point $p\in M$, and it is given by
$\exp_p(\nu_p(Gp))$. 
\item If $\Sigma_1$ and $\Sigma_2$ are two sections of $(G,M)$,
then there exists $g\in G$ such that $g\Sigma_1=\Sigma_2$.
In other words, any two sections differ by an element of the
group.   
\end{enumerate}
\end{prop}

\Pf (a) If $\Sigma$ meets a given orbit $N$ at a point $p$,
then $g\Sigma$ meets $N$ at the point $gp$. This shows that 
$g\Sigma$ meets all the orbits. Moreover, if $g\Sigma$
meets $N$ at a point $q$, then it is perpendicular there,
because $\Sigma$ meets $N$ at $g^{-1}q$ and this is perpendicular
and $G$ acts by isometries. It follows that $g\Sigma$ 
satisfies the two defining conditions of a section. 

(b) Let $\Sigma$ be a section of $(G,M)$.
Given $p\in M$, the orbit $Gp$ meets $\Sigma$ in a point $gp$
for some $g\in G$ by the definition of a section. 
Then $g^{-1}\Sigma$ is a section by (a) and $p\in g^{-1}\Sigma$.

(c)  Let $\Sigma$ be a section. Then 
\[ T_p\Sigma\subset \nu_p(Gp) \]
for every $p\in\Sigma$ by definition of a section. 
Denote by $\Sigma_{reg}=\Sigma\cap M_{reg}$ the open set of regular points of
$M$ that lie in~$\Sigma$.
Since $\dim\nu_p(Gp)$ equals the co-homogeneity of $(G,M)$
for $p\in\Sigma_{reg}$, the above inclusion implies
that $\dim\Sigma$ is not larger than this co-homogeneity.
Recall the submersion $\pi:M_{reg}\to M_{reg}/G$.
Since $\Sigma$ intersects all the orbits,
the restriction $\pi|_{\Sigma_{reg}}:\Sigma_{reg}\to M_{reg}/G$ is surjective. 
It follows that $\dim\Sigma_{reg}\geq \dim M_{reg}/G$.
Since $\dim M_{reg}/G$ is equal to the co-homomogeneity 
of $(G,M)$, we conclude that $\dim\Sigma$
is also equal to this co-homogeneity.

(d) It is clear that the set of regular points in $\Sigma$ is open. 
Suppose, on the contrary, that there exists a non-empty open 
subset $V$ of $\Sigma$ that does not contain regular points
of $(G,M)$. Let $p\in V$ be a point whose isotropy subgroup
$G_p$ has the minimal dimension and the smallest number
of connected components among the points in $V$. 
By the normal slice theorem,
$(G_p)=(G_q)$ for $q\in V$. It follows that $GV\approx Gp\times V$ is a 
submanifold of $M$. If $S$ is the normal slice at $p$, then $T_pS=\nu_p(Gp)$
and $T_p(GV)=T_p(Gp)\oplus T_p\Sigma$. It follows that $GV$ is transversal
to $S$ at $p$. By shrinking $V$, we can assume $S\cap GV$ is a 
submanifold $W$, where $\dim W=\dim\Sigma$. $G_p$ cannot fix all the points 
of $S$ because $p$ is not regular, but it fixes all the points of 
$W$, so the co-homogeneity of $(G_p,S)$ is at least
$\dim W+1=\dim\Sigma+1$. 
The cohomogeneity of a slice is also the co-homogeneity of $(G,M)$, 
which contradicts part~(c).

(e) Let $\Sigma$ be a section.
By part (d), $\Sigma_{reg}$ is dense in $\Sigma$. Thus, by continuity,
it suffices to prove that the second fundamental form 
of $\Sigma$ in $M$ vanishes along $\Sigma_{reg}$.
Let $p\in\Sigma_{reg}$ and consider a normal vector 
$u\in\nu_p\Sigma$. Since $p$ is a regular point, 
$\nu_p\Sigma=T_p(Gp)$, so we can find an element $X$
in the Lie algebra of $G$ such that $X^*_p=u$.
Owing to the polarity of the action, 
$X^*$ is perpendicular to $\Sigma$ everywhere along~$\Sigma$.
Therefore the Weingarten operator $A$ of $\Sigma$
can be written as $\langle A_u v,v\rangle=-\langle\nabla_vX^*,v\rangle=0$
for all $v\in T_p M$, where we have used that $X^*$ is a Killing field.
Hence $A$ vanishes along $\Sigma_{reg}$, as we wanted. 

(f) Let $p\in M$ be a regular point and let $\Sigma$ be a section
through $p$. We have seen that $T_p\Sigma=\nu_p(Gp)$ and
$\Sigma$ is totally geodesic, so $\Sigma\supset\exp_p(T_p\Sigma)
=\exp_p(\nu_p(Gp))$. For the converse inclusion, let $q\in\Sigma$.
By part~(e) and completeness of $\Sigma$, there exists a minimizing geodesic
of $M$, $\gamma:[0,1]\to\Sigma$, with $\gamma(0)=p$ and $\gamma(1)=q$. 
Then $q=\exp_p(\gamma'(0))$ where $\gamma'(0)\in T_p\Sigma$,
proving that $\Sigma\subset\exp_p(T_p\Sigma)=\exp_p(\nu_p(Gp))$. 

(g) Let $p\in\Sigma_1$ be a regular point. There exists 
$g\in G$ such that $gp\in\Sigma_2$. Now 
$g\Sigma_1$ and $\Sigma_2$ are two sections 
through the regular point $gp$, so they 
must coincide by part (f). \EPf 

\medskip

The next result shows that the property 
of being polar is inherited by the slice 
representations. 

\begin{prop}\label{slice-representation-polar}
Let $(G,M)$ be a polar action, and let $p\in M$. 
Then the slice representation at $p$ is also polar. 
In fact, if $\Sigma$ is a section of $(G,M)$ 
through $p$, then $T_p \Sigma$ is a section 
of $(G_p,\nu_p(Gp))$. 
\end{prop}

\Pf Set $K=G_p$ and $V=\nu_p(Gp)$ for convenience. 
The co-homogeneity of the 
slice representation is the same as that 
of the slice action of $K$ on the normal slice $S$ at $p$.  
This shows that $T_p\Sigma$ has the right
dimension of a section of $(K,V)$. 

We claim that $T_p\Sigma$ contains regular
points of~$(K,V)$. In fact, let $\xi$ be a principal
orbit of $G$, and choose a connected component~$\beta$
of~$\xi\cap\Sigma$. Let $\gamma$ be a minimizing geodesic
in $\Sigma$ from $\gamma(0)=p$ to $\gamma(1)\in\beta$.
Then $\gamma'(0)\in T_p\Sigma$ is a regular point
of~$(K,V)$. Next, if we can prove that 
$T_p\Sigma$ is perpendicular to $Kv$ for every $v\in T_p\Sigma$,
this will finish the proof, for
it will follow that, for a $(K,V)$-regular point 
$w\in T_p\Sigma$, 
$T_p\Sigma$ coincides with the normal space of $Kw$
in $V$, and hence $T_p\Sigma$ meets all the 
$K$-orbits in $V$ owing to
Lemma~\ref{exponential-normal-space}.  

So let $v\in T_p\Sigma$. The Lie algebra $\Lk$ 
consists of the elements $X$ of $\Lg$ such that~$X_p^*=0$. We also have that  
$\Lk$ induces Killing fields on $V$ via the 
action of~$(K,V)$; denote them with $()^{**}$. 
The tangent space $T_v Kv$
is spanned by the vectors $X^{**}_v\in V$, where 
$X\in\Lk$. Let $\gamma$ be an integral curve of $v$ through~$p$.
In view of
\begin{align*}
  X^{**}_v&=\frac{\nabla}{dt}\Big|_{t=0}d(\exp(tX))_p(v)\\
          &=\frac{\nabla}{dt}\Big|_{t=0}\frac{d}{ds}\Big|_{s=0}(\exp(tX))(\gamma(s))\\
&=\frac{\nabla}{ds}\Big|_{s=0}\frac{d}{dt}\Big|_{t=0}(\exp(tX))(\gamma(s))\\
&=\frac{\nabla}{ds}\Big|_{s=0}X_{\gamma(s)}\\
          &=\nabla_vX,
            \end{align*}
for $u\in T_p\Sigma$ we have 
\[ \langle X^{**}_v,u \rangle = \langle (\nabla_v X^*)_p, u\rangle
=-\langle A_{X^*_p}v,u\rangle = 0. \]
This shows that $Kv$ is perpendicular to $T_p\Sigma$
and completes the proof. \EPf

\begin{cor}\label{isotropy-subgroup-transitive-sections}
Let $(G,M)$ be a polar action, and let $p\in M$.
Then the isotropy subgroup $G_p$ acts 
transitively on the set of sections of $(G,M)$
through~$p$.
\end{cor}

\Pf Let $\Sigma_1$ and $\Sigma_2$ be two sections
containing the point $p$. According to 
Proposition~\ref{slice-representation-polar},
the slice representation $(G_p,\nu_p(Gp))$ is polar
and $T_p\Sigma_1$, $T_p\Sigma_2$ are two of its
sections. 
By Propostion~\ref{basic-properties-polar-action}(g),
there exists $g\in G_p$ such that $dg_p(T_p\Sigma_1)=T_p\Sigma_2$.
Since $\Sigma_1$ and $\Sigma_2$ are totally geodesic, 
this implies that $g\Sigma_1=\Sigma_2$, as wished.  \EPf

\subsection{The generalized Weyl group}

Let $(G,M)$ be a proper isometric action.
By the normal slice theorem, locally, 
the study of the orbit space near an orbit $Gp$
is reduced to the study of the orbit space
of the action of $G_p$ on the normal slice $S_p$.
Next, we explain how this reduction can be done
globally in the case in which 
$(G,M)$ is a polar action. 

Let $(G,M)$ be a polar action, and let $\Sigma$ be a 
section. 
Denote by $N(\Sigma)$ the normalizer of $\Sigma$ in 
$G$, namely, the subgroup of $G$ consisting of the  
elements that restrict to isometries of $\Sigma$. 
Then the action of $G$ on $M$ restricts to an 
action of $N(\Sigma)$ on $\Sigma$. 
In the following, it will be interesting 
to consider the effectivized action of $N(\Sigma)$
on $\Sigma$; we say an action is \emph{effective}
if the only group element that acts as the identity map
is the identity element in the group. For that purpose, denote by
$Z(\Sigma)$ the centralizer of $\Sigma$ in $G$, namely,
the subgroup of $G$ consisting of the  
elements that restrict to the identity on $\Sigma$. 
Take any regular point $p\in\Sigma$. It is obvious that 
$Z(\Sigma)\subset G_p$, and the reverse inclusion is a 
consequence of the fact that the slice representation at a 
regular point is trivial.  
In particular, $Z(\Sigma)=G_p$ is a closed subgroup of~$G$. 
Note also that $N(\Sigma)$ is a subgroup of the normalizer
of $G_p$ in $G$, $N(\Sigma)\subset N_G(G_p)$. 

The \emph{generalized Weyl group}\index{Weyl group!generalized}
of the polar action $(G,M)$ with respect to the 
section $\Sigma$ is defined to be the quotient group
\[ W(\Sigma)=N(\Sigma)/Z(\Sigma). \]
In Proposition~\ref{weyl-group-properties},
we collect a number of properties
related to the generalized Weyl group.\ft{In the case 
of an isometrically immersed section $\iota:\Sigma\to M$, 
$W(\iota)$ is defined as the normalizer of the image of~$\iota$, quotiented
by its centralizer. One then shows that the $W(\iota)$-action
on $\iota(\Sigma)$ lifts to an action on~$\Sigma$: this is clear on the 
$G$-regular points of $\Sigma$, and follows on the other points by continuity.} 

\begin{rem}\label{horizontal-geodesic-realizes-distance}
Recall that for a proper isometric action~$(G,M)$, where 
$M$ is assumed connected, 
the set of manifold points of $M^*=M/G$ is open and dense, but 
$M^*$ is, in general, not a manifold. Nevertheless,
there is a natural quotient metric space structure on~$M^*$. 
Let $x$, $y\in M^*$. One defines the distance 
$d(x,y)$ to be the distance between the $G$-orbits 
$\pi^{-1}(x)$ and $\pi^{-1}(y)$ in~$M$. Since the $G$-action is proper, 
its orbits are properly embedded submanifolds of $M$, and 
therefore $d(x,y)>0$ for $x\neq y$. It is then clear that
$d$ defines a structure of metric space on~$M^*$. 
Note that $d(x,y)$ is equal to the length of a 
geodesic of $M$ joining a point in $\pi^{-1}(x)$ to a 
point in $\pi^{-1}(y)$, which is horizontal in the sense
that it is orthogonal to every $G$-orbit that it meets;
the initial point of the geodesic
in one of the two orbits can be any chosen point,
by $G$-invariance, but this determines the endpoint in the other orbit.  
\end{rem}

\begin{prop}\label{weyl-group-properties}
Let $(G,M)$ be a polar action admitting a section
$\Sigma$. 
\begin{enumerate}
\item The generalized Weyl group $W(\Sigma)$ is a 
discrete closed subgroup of $N(G_{pr})/G_{pr}$, for some 
principal isotropy group $G_{pr}$ of $(G,M)$.
In particular, $W(\Sigma)$ acts 
properly on $\Sigma$.
\item The inclusion $\iota:\Sigma\to M$ induces a 
map $\bar\iota:\Sigma/W(\Sigma)\to M/G$, which takes a
$N(\Sigma)$-orbit of a point in $\Sigma$ to the 
$G$-orbit of that point. This map is  
is an isometry between the quotient metric spaces.
\item For every $p\in\Sigma$,
$Gp \cap \Sigma = W(\Sigma)p$. 
\item If $\Sigma_1$ and $\Sigma_2$ are sections,
then there exists an isomorphism between
the generalized Weyl groups 
$W(\Sigma_1)$ and $W(\Sigma_2)$ which is uniquely 
defined up to an inner automorphism of $W(\Sigma_1)$.
\end{enumerate}
\end{prop}

\Pf (a) Let $p\in\Sigma$ be a regular point and write $G_p=G_{pr}$. 
Let $S$ be the normal slice at $p$. Then $S$ is an open
neighborhood of $p$ in $\Sigma$. The continuity of the 
action implies that $gp\in S$ for an element 
$g\in N(\Sigma)$ sufficiently close to the identity of 
$N(\Sigma)$. Since $Gp$ is a principal orbit, 
$S$ meets every orbit near $p$
at a unique point, 
so $gp=p$, namely, $g\in G_p=Z(\Sigma)$. 
This argument thus shows that $Z(\Sigma)$ contains an 
open neighborhood of the identity in $N(\Sigma)$, and 
this is equivalent to saying that $Z(\Sigma)$ is an 
open subgroup of $N(\Sigma)$. Hence the quotient
$N(\Sigma)/Z(\Sigma)$ is a discrete Lie group.
Now every discrete subgroup of a Hausdorff topological group
is closed.
The properness of the $W(\Sigma)$-action on $\Sigma$ is immediate from this and 
the properness of the $G$-action on $M$. 

(b) Since $\Sigma$ meets
all the orbits of $G$ in $M$, this map is 
surjective. We claim that the map $\bar\iota$ is also injective.
In order to prove this claim, suppose that $p$, $q\in\Sigma$ 
lie in the same $G$-orbit; we need to prove 
that they lie in the same $N(\Sigma)$-orbit, too. 
We can write $q=gp$ for some $g\in G$.
Then $q$ lies in $g\Sigma$, so $\Sigma$ and $g\Sigma$ are two sections
through the point $q$. 
By Corolllary~\ref{isotropy-subgroup-transitive-sections},
there exists $h\in G_q$ such that $hg\Sigma=\Sigma$. 
It follows that $q=hq=hgp$ where $hg\in N(\Sigma)$,
and this proves the claim. 

We already know that $\bar\iota$ is a continous and bijective 
map. Also, the map $\bar\iota$ is non-expanding (or $1$-Lipschitz), namely,
\[ d(\bar\iota(x),\bar\iota(y))\leq d(x,y) \]
for all $x$, $y\in\Sigma/W$, since every geodesic in $\Sigma$ is a 
geodesic in $M$. 
In view of the denseness
of $G$-regular points in $\Sigma$, 
to finish the proof we need only show that $\bar\iota$ is an 
isometry on $\Sigma_{reg}$. In fact let $x=Np$, $y=Nq$ where 
$p$, $q\in\Sigma_{reg}$. The minimizing geodesic $\gamma$ in $M$ from 
$p$ to $Gq$ is entirely contained in $\Sigma$. Let $r\in\Sigma\cap Gq$
be the endpoint of $\gamma$. Clearly $\gamma$ minimizes the 
distance from $Np$ to $Nr$. Since $Gr=Gq$, by 
the injectivity of $\bar\iota$ we have $Nr=Nq$. Hence 
$d(\bar\iota(x),\bar\iota(y))=\mathrm{Length}(\gamma)=d(x,y)$ as desired.
\mbox{ } \hfill\mbox{ }\EPf

(c) One inclusion is obvious, and the other one is the injectivity of
the map~$\bar\iota$ proved in part~(b).

(d) By Proposition~\ref{basic-properties-polar-action}(g),
there exists an element $g\in G$ such that $g\Sigma_1=\Sigma_2$. 
It is easy to see that $gN(\Sigma_1)g^{-1}=N(\Sigma_2)$
and $gZ(\Sigma_1)g^{-1}=Z(\Sigma_2)$. So the conjugation
by $g$ induces an isomorphism $W(\Sigma_1)\to W(\Sigma_2)$. 
If $g'\in G$ is another element satisfying $g'\Sigma_1=
\Sigma_2$, then $g^{-1}g'\in N(\Sigma_1)$, so 
this element defines an inner automorphism
of $W(\Sigma_1)$ and the 
conjugations by $g$ and $g'$ induce isomorphisms
$W(\Sigma_1)\to W(\Sigma_2)$ that differ by 
that inner automorhism. \hfill\mbox{ }\hfill\EPf

\medskip

\subsection{The orbit space}

Recall that a \emph{Riemannian orbifold} is a length space
locally isometric to the quotient of a Riemannian manifold
by a finite group of isometries (cf.~Lecture~\ref{sec:orbifolds}). For a section
$\Sigma$ of a polar action $(G,M)$,
the generalized Weyl group $W(\Sigma)$ is a discrete
group acting properly on $\Sigma$
(Proposition~\ref{weyl-group-properties}(a)); thus its isotropy subgroups
are finite. It follows from Proposition~\ref{weyl-group-properties}(b)
that the orbit space of a polar action is a Riemannian orbifold. 

Due to Proposition~\ref{weyl-group-properties}(b), the action 
of $(W(\Sigma),\Sigma)$  can also be seen as 
a ``reduction'' of the action~$(G,M)$ to a discrete group action,
namely, one can recover the same orbit space from a much simpler
action of a discrete group action. It easily follows 
from O'Neill's equation~(\ref{oneill}) and Theorem~\ref{thm:hlo}
that a proper and isometric action admits a reduction to a discrete
group if and only if it is a polar action.

Consider for instance the 
case of an orthogonal representation~$(G,V)$.  
In invariant theory, if $(H,W)$ is a reduction of $(G,V)$, that is, 
$W/H$ is isometric to $V/G$, it is a natural question to ask if 
the invariant rings of these representations are isomorphic.
In some special cases this question has an affirmative answer,
namely, polar representations (by Chevalley's theorem) and 
the reduction to the principal isotropy group 
(by Luna-Richardson's theorem~\cite{luna-richardson}).
In~\cite{alexandrino-radeschi} it is shown that the answer 
is also positive in case the isometry preserves the codimension
of the orbits, and in~\cite{alexandrino-lytchak} it is remarked to hold
for infinitesimally polar actions (cf.~section~\ref{sec:orbifolds});
see also~\cite{mendes} for the special case of isometries of $V/G$.
In full generality, the question remains open. 

The reduction principle for orthogonal representations 
was apparently first stated in~\cite{straume}.
In~\cite{gorodski-lytchak-lowcop} a 
systematic study of reductions of orthogonal representations was initiated, 
going much beyond polar representations. As an application
of the results in~\cite{gorodski-kollross-wilking}, there was
shown that an irreducible representation of a compact connected
simple Lie group admits a reduction if and only if it is polar
or it admits a reduction to a finite extension of a torus or a product
of $Sp(1)$-subgroups, and irreducible representations of compact
connected Lie groups admitting reductions to a finite extension of a
torus or a product of $Sp(1)$-subgroups were classified
in~\cite{gorodski-lytchak-toric,gorodski-gozzi}.

\subsection{Examples and classification}

We first discuss the linear case. 
The standard examples of polar representations are the isotropy
representations of symmetric spaces without an Euclidean factor,
sometimes called
\emph{s-representations}. Let $M=G/K$ be a such a symmetric space,
and consider the associated
decomposition on the Lie algebra level $\Lg=\Lk+\Lp$ into the
$\pm1$-eigenspaces of the involution. Then the isotropy representation
of $M$ is equivalent to the adjoint representation of $K$ on $\Lp$. 
Recall that a symmmetric space of compact type and its noncompact
dual have equivalent isotropy representations. By passing to a
covering, we may assume that $K$ is connected, and $M$ is irreducible.
Then the metric is proportional to the Killing form $B$ on~$\Lp$. 
Let $\La$ be a maximal Abelian subalgebra of
$\Lp$ ($\La$ exponentiates to a maximal flat
of $G/K$ through the basepoint). Let $X\in \Lk$, $Y$, $Z\in\La$
be arbitrary. Then
$[X,Y]$ is an arbitrary tangent vector to the $K$-orbit through $Y$ at
the point~$Y$, and 
\[ B([X,Y],Z) = B(X,[Y,Z]) = 0, \]
where we have used the $\ad$-invariance of the Killing form,
and the fact that $\La$ is Abelian. This shows that $\La$ is orthogonal to
$K(Y)$ at $Y$. If $Y$ is a $K$-regular point, then $\dim K(Y)$
equals the codimension of $\La$, so in this case $\La$ coincides
with the normal space $\nu_Y(K(Y))$. It follows from
Lemma~\ref{exponential-normal-space} that $(K,\Lp)$ is a polar
representation. Here the generalized Weyl group coincides with the 
(little) Weyl group of the symmetric space. 

The classification theorem of Dadok~\cite{dadok}
implies that
every polar representation of a compact connected Lie group
is orbit-equivalent to such an s-representation, for some
symmetric space (see also~\cite{kollross2003} for an alternative approach,
and~\cite{eschenburg-heintze2} for a geometric proof in case
of cohomogeneity bigger than two).
The list of representations orbit-equivalent to
an s-representation is not difficult to compile. Besides the cohomogeneity
one case which was previously known, in the irreducible case there is a short
list~\cite{eschenburg-heintze}, and in the general case there is a
description~\cite{bergmann}.

Note that the orthogonal conjugacy of real symmetric matrices to
diagonal matrices proved in basic Linear Algebra courses is
the polarity of the s-representation associated to~$SL(n,\R)/SO(n)$. 

Moving to polar actions on more general Riemannian manifolds,
for a proper and isometric action, a geodesic orthogonal 
to an orbit remains orthogonal to every orbit it meets, due 
to Clairault's lemma. It follows that cohomogeneity one actions 
form a class of polar actions with a flat section.\ft{Either due 
to the brothers Alekseevsky result, or to the broader 
definition of a section, cf.~footnote~1.} It is interesting 
to remark that flat sections of polar actions on symmetric spaces
of compact type equipped with metrics coming from the Cartan-Killing form
were shown to be compact tori (hence properly embedded) 
in~\cite{HPTT}.\ft{A cohomogeneity one
action on a torus can of course have as section a dense irrational torus. 
Even in the case of a cohomogeneity one action on a compact symmetric 
space without flat factor, taking a metric non-proportional 
to the Cartan-Killing form in the reducible case, one can have a
non-embedded section: one such example is the diagonal
action of $SO(3)$ on the product of spheres
$S^2(1)\times S^2(R)$ where $R^2$ is irrational.
(This follows because $\gamma(t)=((\cos t,\sin t,0),(R\sin\frac{t}{R^2},
R\cos\frac{t}{R^2},0))$ is a geodesic normal to the orbits.)} 

A polar action with flat sections is sometimes called \emph{hyperpolar}. 
In~\cite{hermann}, Hermann
 constructed examples of ``variationally complete'' actions
(in the sense of Bott and Samelson, cf.~Lecture~3),
which were later found to be 
hyperpolar, as follows. Let $(G,H)$ and $(G,K)$ be two 
symmetric pairs where, say, $G$ is compact.
Then $G/H$ and $G/K$ are compact symmetric spaces. 
The left-action of $H$ on $G/K$, now called a 
\emph{Hermann action}, is hyperpolar\ft{The left- and right-action of 
$H\times K$ on $G$ is also hyperpolar. These actions generalize
the isotropy action $(K,G/K)$, and the left- and right-action 
$(K\times K,G)$, which had been previously considered by Bott and Samelson.}
(cf.~section~\ref{criterion}).

Kollross classified hyperpolar actions (in particular, cohomogeneity one
actions) on compact irreducible symmetric 
spaces in his thesis, which is published as~\cite{kollross2000}. 
Later in~\cite{kollross2017} he showed that one can remove the 
irreducibility assumption if the cohomogeneity of the action is bigger
than one.\ft{The classification of cohomogeneity one actions on 
reducible symmetric spaces is still open.} His result is that 
an indecomposable\ft{Here an isometric action on a Riemannian manifold
is called \emph{decomposable} if the manifold can be written as a 
Riemannian product, and the action is orbit equivalent to the product 
of isometric actions on the factors, and \emph{indecomposable} otherwise.}
hyperpolar action of cohomogeneity at least two
on a compact symmetric space is orbit-equivalent to a Hermann action.

There are easy examples of polar actions with nonflat sections
on compact symmetric spaces of rank one. A simple example is the
action of the maximal torus $T^n$ of the isometry group $SU(n+1)$
of complex projective space $\C P^n$, which is polar with sections isometric
to a totally geodesic $\R P^n$. According to the classification
result of Podest\`a and Thorbergsson~\cite{podesta-thorbergsson}, the polar
actions on classical symmetric spaces of rank one are induced from 
certain polar actions on spheres using the Hopf fibration.
In the case of the Cayley projective plane, there is no Hopf fibration,
and their analysis found four polar actions with cohomogeneity one, and
four polar actions with cohomogeneity two; a further
polar action of cohomogeneity two which
was overlooked in~\cite{podesta-thorbergsson} was later
found to be polar in~\cite{gorodski-kollross} (it is given by a maximal
subgroup of the isometry group of $\Ca P^2$ whose Lie algebra
is not regular). 

No examples of polar actions with nonflat sections on irreducible compact
symmetric spaces of rank bigger than one were known,
so eventually the question of their existence became a folklore problem.
Many special cases were examined, most notably in the case of a
Hermitian symmetric space in~\cite{biliotti}. Finally, Lytchak and 
Kollross proved that they do not exist in~\cite{kollross-lytchak}. 

One direction in which to extend the above results is 
to pass from compact symmetric spaces to compact 
non-negatively curved manifolds. Using the theory of Tits buildings, 
Fang, Grove and Thorbergsson proved in~\cite{fang-grove-thorbergsson} 
that a polar action of a compact Lie group on a simply connected
compact positively curved manifold of cohomogeneity at least two is
equivariantly
diffeomorphic to a polar action on a compact rank-one symmetric space.

Another direction to go is to consider polar actions on 
symmetric spaces of non-compact type. Here much fewer results are known.
On real hyperbolic space $\R H^n$ the classification
is the work of Wu~\cite{wu}.
On complex hyperbolic space $\C H^n$ it has been achieved
by D\'\i az-Ramos, Dom\'\i nguez-V\'azques and Kollross~\cite{DDK}
(previously, the case of $\C H^2$ was dealt with in~\cite{BD}).
Next, there is a general description and a partial classification
of cohomogeneity one actions on irreducible symmetric spaces
of non-compact type, due to Berndt and Tamaru~\cite{BT}. They fall into
two cases, namely, either there is a unique singular orbit and the other orbits are distance tubes around the singular orbit, or there are no singular orbits
and the orbits define a regular Riemannian foliation. More generally,
hyperpolar actions with no singular orbits on symmetric spaces
of non-compact type are classified in~\cite{BDT}. Regarding
polar actions on higher rank symmetric spaces of non-compact
type, there is a classification of those polar actions with a
fixed point~\cite{DK}. In another work, Kollross~\cite{kollross-duality}
used Cartan duality between
symmetric spaces of compact type and non-compact type, which maps
an action of reductive algebraic group to a dual action, and showed
that the dual action shares many properties in common with
the original action, for instance (hyper)polarity; in this way, he obtained
a number of new results on polar and hyperpolar actions on
noncompact symmetric spaces. See the survey~\cite{DDS} for more details 
on polar actions on symmetric spaces of non-compact type. 

All examples of polar actions discussed above turn out to be on
symmetric spaces. In~\cite{grove-ziller} an algorithm to recover a
polar action, up to equivariant diffeomorphism,
from the specification of the orbit space and the isotropy
groups along the strata, which must satisfy certain compatibility conditions,
is given, which allows to construct polar actions on a variety of non-symmetric
manifolds. This algorithm was used in~\cite{gozzi} to classify polar
actions on compact simply-connected manifolds up to dimension~$5$,
up to equivariant diffeomorphism (see also~\cite{mucha} for examples
of polar actions on non-symmetric manifolds).

\subsection{Polarity and the 
integrability of the distribution of normal spaces
to the principal orbits}

For a polar action $(G,M)$, the sections are clearly integral 
manifolds of the tangent distribution $\Hd$ on the regular set $M_{reg}$,
defined by $\Hd_p = \nu_p(Gp)$. Note that this is the
horizontal distribution of the Riemannian submersion
$M_{reg}\to M_{reg}/G$.  
As early as 1987, Palais and Terng
conjectured that if $\Hd$ is integrable then its integral manifolds
can be extended to sections,
see~\cite[Remark~3.3]{palais-terng}, where they call
$\Hd$ the \emph{principal horizontal distribution}. 
 As of now, the 
conjecture has been verified and there are different proofs 
in the literature, using different techniques, 
and with differents degrees of generality. We shall now sketch 
the approach of Heintze, Liu and Olmos~\cite[Appendix~A]{heintze-liu-olmos}.

For ease of presentation, we start with the linear case.

\begin{prop}
Let $(G,V)$ be an orthogonal representation of a compact Lie group
$G$ on an Euclidean vector space $V$. If the principal
horizontal distribution $\mathcal H$ is integrable, then~$(G,V)$ 
is polar. 
\end{prop}

\Pf Let $L$ be a leaf of $\Hd$. 
It follows from the argument in
Proposition~\ref{basic-properties-polar-action}(e) that $L$ 
is totally geodesic. Therefore it is a non-empty open
subset of an affine subspace $\Sigma$ of $V$. We claim that $\Sigma$ 
is a section of $(G,V)$. In fact, $\Sigma=T_pL=\nu_p(Gp)$ 
for $Gp\in L$, so $\Sigma$ meets all $G$-orbits by 
Lemma~\ref{exponential-normal-space}. To see that it meets 
always orthogonally, let $X\in\Lg$, and let $\gamma$ be
any horizontal geodesic with $\gamma(0)=p\in L$. Then the 
image of $\gamma$ is contained in~$\Sigma$  and $J:=X^*\circ\gamma$
is a Jacobi field along~$\gamma$. Since $\Sigma$ is totally geodesic,
also the horizontal component $J^h$ of $J$ is a 
Jacobi field. Now $J^h$ vanishes on a neighborhood of
$t=0$, so $J^h$ vanishes identically. This shows that 
$X^*_{\gamma(t)}$ is orthogonal to~$\Sigma$ for all~$t\in\gamma^{-1}(V_{reg})$. 
Since $\gamma$ is arbitrary, $\Sigma$ is orthogonal to all 
principal orbits it meets, and hence to all orbits it meets, by 
denseness. \EPf

\medskip

In the case of an arbitrary complete Riemannian manifold~$M$,
we can start the proof with the same argument, but the main issue
is the completeness of the leaf $L$ of $\Hd$. For that, we 
shall use Hermann's extension of Cartan's criterion for the
existence of a totally geodesic submanifolds with given tangent space
at one point (cf.~section~\ref{cartan-hermann}).

\begin{thm}\label{thm:hlo}
Let $(G,M)$ be a proper isometric action, where $M$ is connected and complete.
Assume that the principal horizontal distribution $\Hd$ is integrable. 
Then $(G,M)$ is polar in the broader sense (cf.~footnote~1). 
\end{thm}

\Pf Fix a $G$-regular point $p\in M$. Let $L$ be maximal leaf of $\Hd$
through~$p$, and put $S=T_pL$. Consider a $S$-admissible 
once-broken geodesic $\gamma:[0,\ell]\to M$ emanating from $p$
(cf.~section~\ref{cartan-hermann} for this terminology) 
that does not meet the set of singular points, that is,
$\gamma(t)$ lies either in a principal or in an exceptional $G$-orbit,
for all $t\in[0,\ell]$. Then $\Hd$ is defined along $\gamma$, and it is 
an auto-parallel distribution; in particular, $\Hd_{\gamma(t)}$ is invariant 
under $R_{\gamma(t)}$ for all~$t$. It follows that~(\ref{condition})
is satisfied along $\gamma$. Next, we consider a 
$S$-admissible once-broken geodesic emanating from $p$ which meets
the singular set at finitely many points. It follows from
 Lemma~\ref{extension} below that~(\ref{condition})
is satisfied along $\gamma$. Finally, we observe that the classes
above are dense in the set of $S$-admissible once-broken geodesic 
emanating from $p$. We deduce that~(\ref{condition}) is satisfied 
in general, by continuity. Now Theorem~\ref{hermann} yields a 
complete totally geodesic isometric immersion $\Sigma\to M$ such that 
$T_p\Sigma=S$, which is clearly a section. \EPf

\medskip

For $v\in TM$, write $\gamma_v(t)=\exp(tv)$. 

\begin{lem}\label{extension}
Assume $\Hd$ is integrable, and let $q\in M$ be a singular point.
Then the slice representation at~$q$ is polar and, for any section
$\Sigma_0$ of it, and any $G_q$-regular $v\in \Sigma_0$, $\Sigma_0$ is the 
parallel transport of $\Hd_{\gamma_v(t)}$ along $\gamma_v$
from $\gamma_v(t)$ to $q$, for small $t>0$.
\end{lem}

\Pf We first prove that the slice representation at~$q$ is polar. 
Let $\gamma_v$ be any geodesic such that $v\in\nu_q(Gq)$ is a 
regular point for the slice representation. Then $\gamma_v(t)$ 
is $G$-regular for small $t>0$. Take a sequence $t_n\ssearrow0$;
by compactness, we may pass to a subsequence and assume
that $\Hd_{\gamma_v(t_n)}$ converges to a subspace $\Sigma_0\in\nu_q(Gq)$. 
Note that $v\in\Sigma_0$. Now  $\langle \nabla_uX^*|_{\gamma_v(t_n)},w\rangle = -\langle A_{X^*_{\gamma_v(t_n)}}u,w\rangle=0$, for all $X\in\Lg$, $u$, $w\in\Hd_{\gamma_v(t_n)}$. 
By continuity,  $\langle X^{**}_u,w\rangle = \langle \nabla_uX^*|_q,w\rangle=0$
for all $X\in\Lg_q$, $u$, $w\in\Sigma_0$. It follows that $\Sigma_0$ is a 
section for the slice representation. 

Next, recall that $\exp_q:\nu_q^{\leq\epsilon}(Gq)\to S_q$ is a
$G_q$-equivariant diffeomorphism, 
where $S_q$ is the normal slice and $\epsilon>0$ is small. 
It follows that 
\begin{equation}\label{exp-tangent-orbits}
d(\exp_q)_{tv}(T_vG_q(v))=T_{\gamma_v(t)}(G_q\gamma_v(t))
\end{equation} 
for $t>0$ small. Put $\Sigma=\exp_q(\Sigma_0)$. 
Since the metric $(\exp_q^*g)_{tv}\to g_q$ uniformly 
on compact subsets as $t\to0$,  taking orthogonal complements 
in~(\ref{exp-tangent-orbits}) we obtain that the distance
in $\Lambda^c(T_{\gamma_v(t)}S_q)$ satisfies
$d(T_{\gamma_v(t)}\Sigma,\Hd_{\gamma_v(t)})\to0$ as $t\to0$,
where $c$ is the cohomogeneity of~$(G,M)$. Together with 
$T_{\gamma_v(t)}\Sigma\to\Sigma_0$,  this implies that 
$\Hd_{\gamma_v(t)}\to\Sigma_0$, and hence $\Sigma_0$ is the parallel
transport of~$\Hd_{\gamma_v(t)}$, by continuity. \EPf

\section{Lecture 2: Orbifold points}\label{sec:orbifolds}

Historically, orbifolds first arose as manifolds
with singular points long before they were formally defined.
Definitions of orbifolds were given by Satake in 1950s
in the context of automorphic forms (``$V$-manifolds''),
by Thurston in the 1970's in the context of $3$-manifolds
(when together with students he coined the name orbifold;
the notions of orbifold coverings and orbifold fundamental groups
are also due to him), and by Haefliger in he 1980's in the context
of CAT($\kappa$)-spaces (``orbihedrons'').
Formally speaking, today there are two
accepted ways to define orbifolds: by means of orbifold atlases,
and this can be done in the topological, differentiable or
Riemannian category; or by means of Lie grupoids, which yields a
slightly more general definition, albeit less geometrical.
A more direct definition, in the restricted Riemannian setting, has
been proposed by Lytchak (see~\cite{lytchak-thorbergsson}).
Herein we follow this approach and sketch some of the main
ideas. A fuller account can be found in~\cite{lange}. 

\subsection{Riemannian orbifolds}

An intrinsic metric space
$X$ is called a \emph{Riemannian orbifold} if every 
point $x\in X$ admits a neighborhood $U$ isometric to a 
quotient $M/\Gamma$, where $M$ is a Riemannian 
manifold and $\Gamma$ is a finite group of isometries.
In this definition, $M$ is endowed with the induced 
intrinsic metric, and $M/G$ with the corresponding
quotient metric, which measures the distance between orbits in $M$. 
Recall that a \emph{intrinsic metric space}
is a special kind of metric space in which the
distance between any pair of points can be realized as the infimum of the 
lengths of all rectifiable paths connecting these points. 

A Riemannian orbifold
is called \emph{good} if it is globally isometric to the
quotient space of a Riemannian manifold
by a discrete group of isometries, and \emph{bad} otherwise.

Lemma~\ref{local-rep} below shows that a Riemannian orbifold
$X$ is locally represented as a quotient in a unique way. 
We first recall a special case proved in~\cite{swartz}.
Let $V$ be an Euclidean space and denote by $S(V)$ its unit sphere. 
If $\Gamma$, $\Gamma'$ are two finite subgroups of $O(V)$ which are
conjugate, then the orbit spaces $V/\Gamma$, 
$V/\Gamma'$ are isometric.
In fact, if $\Gamma'=f\Gamma f^{-1}$ for some $f\in O(V)$, then 
there is an induced isometry
\begin{diagram}
V & \rTo^f & V \\
\dTo^\pi && \dTo_{\pi'} \\
V/\Gamma& \rDotsto_{\bar f} & V/ \Gamma'
\end{diagram}
given by $\bar f(\Gamma v)=\Gamma' f(v)$. 
Conversely:

\begin{lem}\label{lin}
If $V/\Gamma$, 
$V/\Gamma'$ are isometric then $\Gamma$, $\Gamma'$ are 
conjugate in $O(V)$.
\end{lem}

\Pf We proceed by induction on $n=\dim V$.
In the initial case of $n=1$, $V\cong\mathbb R$ and
the only possibilities for $\Gamma$ are $\{1\}$ and 
$\{\pm1\}$, which yield $\mathbb R$ and $[0,+\infty)$, resp., 
non-isometric orbit spaces, so we are done. 

Assume now $n\geq2$.  
It is enough to work with $X=S(V)/\Gamma$,
$X'=S(V)/\Gamma'$. Suppose $F:X\to X'$ is an isometry.  
Let $x\in X$ be such that $\Gamma\cdot x$ and $\Gamma'\cdot F(x)$ are 
principal orbits. Choose points $p\in \pi^{-1}(x)$, $p'\in\pi'^{-1}(x')$
and open neighborhoods $U_p$, $U_{p'}$, $U_x$, $U_{x'}$ such that 
$\pi|_{U_p}:U_p\to U_x$, $\pi'|_{U_p'}:U_{p'}\to U_{x'}$ are isometries
and $F(U_x)=U_{x'}$. Then $(\pi'|_{U_{p'}})^{-1}F\pi:U_p\to U_{p'}$ is an 
isometry, where we view $\pi:S(V)\to X$, $\pi':S(V)\to X'$; 
since $S(V)$ is a space of constant curvature, we  
can (uniquely) extend it to a global isometry $f:S(V)\to S(V)$. 
Let $\bar f:X''\to X'$ be the isometry induced 
on the level of quotients, where 
$\Gamma'':=f^{-1}\Gamma'f$
and $X''=S(V)/\Gamma''$. 
Then $\pi''=\bar f^{-1}\pi'f:S(V)\to X''$. 
Therefore, identifying $X\cong X''$ using the isometry $\bar f^{-1}F$, 
we get $\pi''|_{U_p}=\pi|_{U_p}$.
We will show that $\Gamma''=\Gamma$ as subgroups of $O(V)$. 

It suffices to prove that:
\begin{enumerate}
\item[(a)] $\pi$ is completely determined by its restriction to an
open neighborhood $U_p$ of a $\Gamma$-regular point~$p$.
\item[(b)] $\Gamma$ is completely determined by~$\pi$.
\end{enumerate}
Since $\pi$ is a local isometry on the principal stratum, it is determined
along any unit speed geodesic $\gamma$ in $S(V)$ emanating from $p$,
until $\gamma$ reaches a non-regular point, say $q=\gamma(t_0)$ for 
some $t_0>0$. Now $\dot\gamma(t_0)$ belongs to the unit sphere $S_q$ of
$T_qS(V)$. The space of directions $\Sigma_yX$ for $y=\pi(q)\in X$ is 
isometric to $S_q/\Gamma_q$. Since $\dim T_qS(V)=n-1$, the action of $\Gamma_q$
on $S_q$ is known by the induction hypothesis. It follows that the exit
direction of $\pi\circ\gamma$ from~$y$ is known and thus $\pi$ is 
determined along $\gamma$ beyond $t_0$; this proves~(a). Finally,
the elements of~$\Gamma$ are in bijective correspondence with the 
points in $\pi^{-1}(x)$ via the map $\gamma\mapsto\gamma(p)$. 
For each $\gamma\in\Gamma$, we have a commutative 
diagram:
\begin{diagram}
U_p\\
\dTo^\pi & \rdTo^\gamma\\
U_x & \lTo_\pi &\gamma(U_p)=U_{\gamma(p)}
\end{diagram}
Since $\gamma$ is an isometry of $S(V)$, using~(a)
this completely determines it 
as an element of $O(V)$. Hence~(b) is proved. \EPf

\begin{lem}\label{local-rep}
Every isometry $F:M/\Gamma\to M'/\Gamma'$ is 
\emph{locally} induced by a locally defined isometry $f:M\to M'$. Namely, for 
every $x\in M/\Gamma$, there exist connected open 
neighborhoods $U$, $U'$ of $x$, $x'=F(x)$ of the form 
$V/\Gamma_p$, $V'/\Gamma'_{p'}$, where 
$V$, $V'$ are normal neighborhoods of $p\in\pi^{-1}(x)$,
$p'\in\pi^{-1}(x')$, resp., and $U'=F(U)$ ($\pi:V\to V/\Gamma_p$,
$\pi':V'\to V'/\Gamma'_{p'}$ denote the canonical projections).
Moreover 
$F\circ \pi = \pi'\circ f$ for some isometry $f:V\to V'$ with $f(p)=p'$
and $\Gamma'_{p'}=f\Gamma_pf^{-1}$.
\end{lem}

\Pf Let $G=\Gamma_p$, $G'=\Gamma'_{p'}$. 
By restriction to slices we have an isometry $F:V/G\to V'/G'$
with $F(x)=x'$. Here $V$, $V'$ can be taken to be metric balls 
of the same, small radius, around $p$, $p'$, resp. Consider the actions
of $G$, $G'$ on $T_pV$, $T_{p'}V'$, resp. Then
there is an isometry
\[ T_pV/G\cong T_x(V/G)\to 
T_{x'}(V'/G')\cong  T_{p'}V'/G', \]
which we denote by $dF_x$. By Lemma~\ref{lin}, there is an 
isometry $\varphi:T_pV\to T_{p'}V'$ such that 
\begin{diagram}
T_pV & \rTo^\varphi & T_{p'}V' \\
\dTo^{d\pi_p}&&\dTo_{d\pi'_{p'}} \\
T_pV/G & \rTo_{dF_x} & T_{p'}V'/G'
\end{diagram}
is commutative. Since the Riemannian 
exponential maps $\exp_p:T_pV\to V$, $\exp_{p'}:T_{p'}V'\to V'$
are $G$-, $G'$-equivariant diffeomorphisms, resp., we can define 
an equivariant diffeomorphism $f:V\to V'$ by 
\begin{diagram}
T_pV & \rTo^\varphi & T_{p'}V' \\
\dTo^{\exp_p}&&\dTo_{\exp_{p'}} \\
V & \rDotsto_f & V'
\end{diagram}
Finally, there is an induced map
\begin{diagram}
V & \rTo^f & V' \\
\dTo^{\pi}&&\dTo_{\pi'} \\
V/G& \rDotsto_{\bar f} & V'/G'
\end{diagram}

We claim that $\bar f=F$. In fact, for a geodesic 
$\gamma(t)=\exp_pt\dot\gamma(0)$ in~$V$:
\begin{align}\label{barf} 
\bar f\pi\gamma(t) & = \bar f \pi \exp_p t \dot\gamma(0) \\ \nonumber
    &=\pi'\exp_{p'} t \varphi(\dot\gamma(0))
\end{align}
by commutativity of the last two diagrams. 
Since $x$ is a fixed point of $G$, 
$\pi\circ\gamma$ is a metric geodesic of $V/G$.
Thus it is mapped under $F$ to a metric geodesic 
emanating from~$x'$, namely, $\pi'\circ\gamma'$, 
where $\dot\gamma'(0)=\varphi(\dot\gamma(0))$:
\begin{equation}\label{F}
 F\pi\gamma(t) =\pi'\gamma'(t)=\pi'\exp_{p'}t
\varphi(\dot\gamma(0)). 
\end{equation}
Comparison of~(\ref{barf}) and~(\ref{F}) proves the claim. 

It follows $f:V\to V'$ is a local isometry on the regular set and thus,
by continuity, an isometry everywhere. Now the groups $G'$,
$fGf^{-1}$ acting on $V'$ are orbit-equivalent. If not coincident,
they generate a strictly larger group with the same (finite) orbits 
and thus non-trivial principal isotropy groups, a contradiction
(since the slice representation at regular points must be trivial). 
Hence $G'=fGf^{-1}$. \EPf

\medskip

Let $X$ be a Riemannian orbifold and let $x\in X$. 
Locally represent $X$ around $x$ as a quotient $M/\Gamma$
and write $x=\Gamma p$ for some $p\in M$. Since the isotropy group
$\Gamma_p$ acts by isometries on $M$, it can be viewed as 
a subgroup of the orthogonal group $O(T_pM)$.  
Moreover, it follows from Lemma~\ref{local-rep} that 
the congruence class of $\Gamma_p$ is independent of the local representation
of $X$ as a quotient. After identification $T_pM\cong\mathbb R^n$, 
we get a congruence class of subgroups of $O(n)$, called 
the \emph{local group} of $X$ at $x$ and denoted by $\mathrm{Iso}_x(X)$. 
A point $x\in X$ is called a \emph{manifold point} of $X$ 
if $\mathrm{Iso}_x(X)=\{1\}$. 

\begin{eg}
Let a Lie group $G$ act by isometries on a Riemannian manifold.
Then the orbit space has a canonical structure of Riemannian 
orbifold in the following two cases:
\begin{enumerate}
\item[(a)] $G$ is discrete and the action is proper (such orbifolds
are called \emph{good} or \emph{developable}; non-good orbifolds
are also called~\emph{bad});
\item[(b)] $G$ is compact and connected and all orbits 
have the same dimension.
\end{enumerate}
\end{eg}
 
An \emph{orbi-covering} is a continuous map $\pi:Y\to X$ between 
Riemannian orbifolds where every $x\in X$ admits a neighborhood
$U\cong \tilde U/G$ such that every component~$V$ of $\pi^{-1}(U)$
must be of the form $\tilde V/H$, with $H\subset G$,
and $\pi|_V:V\to U$ lifts to a $H$-equivariant 
homeomorphism $\tilde V\to\tilde U$.
Here $G$ is the local group $\mathrm{Iso}_x(X)$ of $x$ and $H$ is the
local group
$\mathrm{Iso}_y(Y)$ at a point $y\in\pi^{-1}(x)$. 
If we identify $\tilde V\cong\tilde U$ 
via this map, the local representation of the covering map 
is $\tilde U/H\to\tilde U/G$ for $H$ a subgroup of $G$.

 It is a fact
that every connected orbifold $X$ admits a \emph{universal 
orbi-covering} $\tilde X$, unique up to equivalence, with the property
that it orbi-covers any other orbi-covering space of $X$. 
The \emph{orbifold fundamental group} of $X$ is the group
$\pi_1^{orb}(X)$ of deck transformations of the universal 
orbi-covering; it acts simply transitively on the fibers
of this orbi-covering. In general, one can write a presentation of the
orbifold fundamental group
in terms of its usual fundamental group and its strata of codimension~$1$
and~$2$ (cf.~\cite{davis}).
The orbifold fundamental group is a refinement of the 
usual fundamental group in the sense that 
an orbifold can be simply-connected
in the topological sense without being simply-connected
in the orbifold sense.

\begin{eg}
Let the cyclic group $\mathbb Z_m$ act by rotations around a fixed axis
on the sphere $S^2$. The orbit space $X$ is a Riemannian
orbifold and topologically
a $2$-sphere but $\pi_1^{orb}(X)\cong\mathbb Z_m$. 
\end{eg}

\begin{eg}
There is a Riemannian orbifold structure $X_{m,n}$ on $S^2$ with 
exactly two non-manifold (conical) points whose local groups
are respectively $\mathbb Z_m$ and $\mathbb Z_n$. Then 
$\pi_1^{orb}(X_{m,n})=\mathbb Z_d$ where $d$ is the greatest common divisor of
$m$ and $n$. The orbifold $X_{m,n}$ is good if and only if $m=n$. 
In particular, $X_{m,1}$ for $m>1$ is called a \emph{teardrop},
and the bad $2$-orbifold depicted in Fig.~1 is the quotient
of $X_{2,1}$ by a reflection. 
\end{eg}

\begin{center}
\begin{tabular}{c}
\includegraphics[width=1.5in]{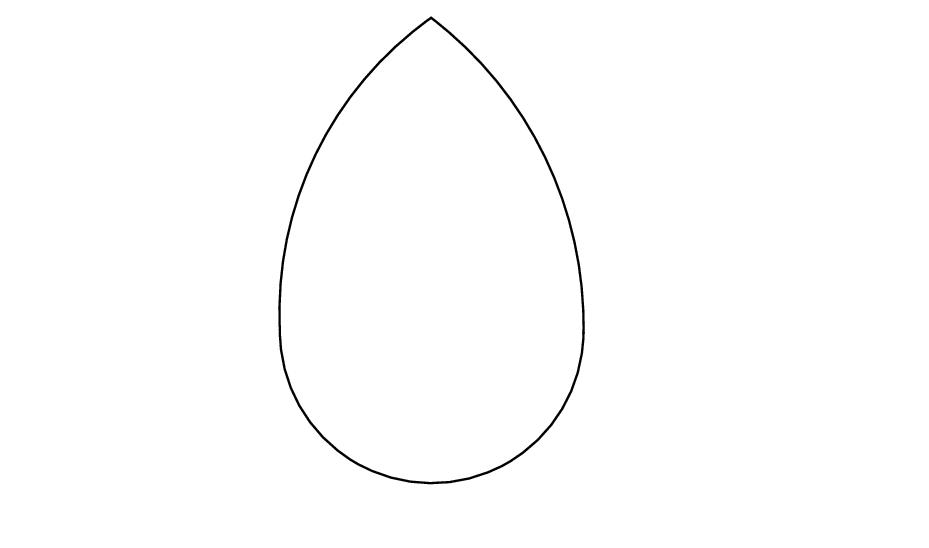}\\
\textsc{Fig.~1:} A bad $2$-orbifold
\end{tabular}
\end{center}

\subsection{Orbifold points in orbit spaces}

Consider a proper and isometric action  of a Lie group $G$ on a 
Riemannian manifold $M$, and let $X=M/G$ be its 
orbit space. The set $X_{reg}$ of regular points of $X$ is exactly the 
set of points that have neighborhoods isometric to Riemannian manifolds,
whereas the slightly larger set $X_{orb}$ of orbifold points of $X$
consists of the set of points that have neighborhoods isometric to  
quotients of Riemannian manifolds by finite groups of isometries.
The main goal of this lecture is to prove the following 
result~\cite{lytchak-thorbergsson}:

\begin{thm}[Lytchak and Thorbergsson]\label{lt1}
Let $(G,M)$ be a proper isometric action of a Lie group 
$G$ on a complete Riemannian manifold $M$, and consider 
its orbit space $X=M/G$. A point $x=Gp\in X$ is an orbifold point 
if and only if the slice representation at~$p\in M$ is
polar.
\end{thm}

It follows from this theorem that $X_{orb}$ contains 
all strata of codimension at most two, since representations
of cohomogeneity at most two are always polar~\cite{hsiang-lawson}.  

\begin{lem}\label{h-lambda}
Let $(M,g)$ be a Riemannian manifold, fix a point $p\in M$,
and consider the locally defined family
of homotheties~$\{h_\lambda\}_{0\leq\lambda\leq1}$, given by 
\[ h_\lambda(\exp_pv)=\exp_p(\lambda v) \]
for $v\in T_pM$. Then the Riemannian 
metrics $g_\lambda:=\frac1{\lambda^2}h_\lambda^*g$ converge 
smoothly on compact neighborhoods of~$p$ to 
$g_0:=(\exp_p^{-1})^*g_p$ as $\lambda\to0$.
\end{lem}

\Pf The proof reduces to a standard calculation in normal coordinates. Write 
$h_\lambda=\exp_p\circ \tilde h_\lambda\circ \exp_p^{-1}$, where 
$\tilde h_\lambda(v)=\lambda v$. It is equivalent to show that
the metrics $\frac1{\lambda^2}\tilde h_\lambda^*\exp_p^*g=
(\exp_p^*g)\circ\tilde h_\lambda$
smoothly converge to the flat metric $g_p$ on 
compact neighborhoods of $0$ in $T_pM$.

For $v\in T_pM$ and $e_1,\ldots,e_n$ an orthonormal basis of
$(T_pM,g_p)$, we have
$(\exp_p^*g)_{\lambda v}(e_i,e_j)=g_{ij}(\lambda v)=:G_{ij}(\lambda,v)$,
where $g_{ij}$ are the coefficients of the metric in the induced
normal coordinates.
But on a compact neighborhood $K$ of $0$ and all $i$, $j$, the smooth function
$G_{ij}$ and all its partial derivatives in~$v$ converge uniformly
to their value at $\lambda=0$ as $\lambda\to0$, that is,  
$g_{ij}\circ\tilde h_\lambda$ converges in the $C^\infty$-topology on~$K$
to the constant function $\delta_{ij}$. \EPf

\medskip

\textit{Proof of Theorem~\ref{lt1}.} Let $S$ be a normal slice at~$p$, and
note that the homotheties $h_\lambda$ restrict to $S$ and induce isometries
$h_\lambda:(S,g_\lambda)\to(h_\lambda(S),\frac1{\lambda^2}g)$,
where we have 
used the notation of Lemma~\ref{h-lambda}. The $G$-action preserves the
metrics $g_\lambda$. For a regular point $q\in S$, denote by
$\bar\kappa_{g_\lambda}(q)$ the supremum of all sectional curvatures of the metric
induced by $g_\lambda$ on the local quotient $GS/G$ at $Gq$. Owing to
Lemma~\ref{h-lambda}, $\bar\kappa_{g_\lambda}(q)\to\bar\kappa_{g_0}(q)$ as $\lambda\to0$. On the other hand,
\begin{align*}
  \bar\kappa_{g_\lambda}(q)  &= \bar\kappa_{\frac1{\lambda^2}g}(h_\lambda(q))\\
  &=\lambda^2\bar\kappa_g(h_\lambda(q))\\
  &=\frac1{||v||^2}d(h_\lambda(q),p)^2\bar\kappa_g(h_\lambda(q)), 
\end{align*}
where $q=\exp_pv$ for $v\in\nu_p(Gp)$, so $d(h_\lambda(q),p)=\lambda d(q,p)=\lambda ||v||$. 

If $x$ is an orbifold point, the sectional curvatures of 
$M_{reg}/G$ near $x$ are locally bounded, so 
\[ \bar\kappa_{g_0}(q)=\frac1{||v||^2}\lim_{\lambda\to0} d(h_\lambda(q),p)^2\bar\kappa_g(h_\lambda(q))=0. \]
It follows that the orbit space of the slice representation at~$p$ 
is flat at regular points. By O'Neill's formula~(\ref{oneill}), the principal 
horizontal distribution of the slice representation is integrable, and hence 
the slice representation is polar. 

Conversely, assume that the slice representation at~$p$ is polar, and let 
$\Sigma$ be a section with associated Weyl group~$W$. Let $N=\exp_p(\Sigma)\cap S$,
where $S$ is a normal slice at $p$. Then $W$ acts on $N$, and we shall define 
a~$W$-invariant Riemannian metric $\tilde g$ on $N$ such that $N/W$ 
is isometric to a neighborhood of $x$ in $X$. 

For $q\in S$, put $\mathcal V_q=(T_qN)^{\perp g_0}$ and 
$\mathcal H_q=(\mathcal V_q)^{\perp g}$. These are smooth distributions,
with $\mathcal V_q\supset T_q(Gq)$ (since $\Sigma$ is $g_p$-orthogonal
to $T_v(G_pv)$ for $v\in\Sigma$), $\mathcal H_q\subset T_q(Gq)^{\perp g}$,
and the latter inclusion is an equality if $q$ is a regular point.  
Let $P_q$ denote the orthogonal projection of $T_qM$ onto $\mathcal H_q$, and define
\[ \tilde g_q(u,v)=g(P_q(u),P_q(v)) \]
for $u$, $v\in T_qN$. Since $P_q:(T_qN,\tilde g_q)\to(\Hd_q,g)$ is an isometry, 
the projection $\pi:(N,\tilde g)\to M/G$ preserves the lengths of all curves 
contained in the regular set of~$(G,M)$. Since $\pi$ is $W$-invariant, 
the action of $W$ on $N$ preserves the length of curves in $N\cap M_{reg}$;
by continuity, $W$ acts by isometries on $N$. It follows that 
$N/W\to M/G$ is an isometric embedding onto a neighborhood of~$x$,
where $N/W$ is a Riemannian orbifold. \EPf

\medskip

\subsection{Applications}

A Riemannian orbifold comes along with a canonical stratification
given by the connected components of the set of points with the
same local group. Each stratum is a connected Riemannian manifold,
which is locally convex with respect to the ambient metric.
The closure of any stratum is a union of
strata. Any Riemannian orbifold $B$ can be written as a quotient of a
Riemannian manifold (the orthonormal frame bundle
of $B$) by an almost free isometric action of a compact Lie group.
The canonical stratification of $B$ is then exactly the stratification
by orbit type.

Let $X=M/G$ be the orbit space of a proper and isometric
action~$(G,M)$. 
The boundary $\partial X$ of~$X$ (in the sense of Alexandrov
geometry) is the closure of the union of strata
of codimension $1$. Since $X_{orb}$ contains all strata
of $X$ of codimension at most~$2$, $X_{orb}$ has non-empty
boundary if and only if $X$ has non-empty boundary.
A Riemannian orbifold with a non-empty boundary can be doubled.
It follows that a Riemannian orbifold $B$ has
$\partial B\neq\varnothing$ if and only if $\pi_1^{orb}(B)$
contains a reflection. 

\begin{eg}
Let $X$ be the quotient of $S^2$ by the reflection across the equator.
Then $X_{orb}=X$ and $\pi_1^{orb}(X)\cong\mathbb Z_2$.
\end{eg}

\begin{rem}\label{rem1}
  Let $(G,M)$ be a proper and isometric action, where $G$ is connected. 
Denote by $B_0$ the subset of points in~$X=M/G$ representing non-singular
$G$-orbits (that is, principal
and exceptional $G$-orbits), so that $B_0=M_0/G$, where
$M_0$ is the union of principal and exceptional $G$-orbits in~$M$.
Then $B_0$ has the structure of a Riemannian orbifold. 
We lift~$(G,M)$ to a proper and isometric action of the simply-connected
covering Lie group $\tilde G$ of $G$ on the simply-connected
Riemannian covering $\tilde M$ of $M$~\cite[Thm.~I.9.1]{bredon}. Then $B_0=\tilde M_0/\tilde G$.
Since all $\tilde G$-orbits have the same dimension in $\tilde M_0$,
there is an epimorphism $\pi_1(\tilde M_0)\to\pi_1^{orb}(B_0)$~\cite{salem}.
Since the union of singular orbits in $\tilde M_0$ has
codimension at most~$2$~\cite[Thm.~IV.3.8]{bredon}, 
$\pi_1(\tilde M_0)=\pi_1(\tilde M)=\{1\}$,
so also $\pi_1^{orb}(B_0)=\{1\}$.
\end{rem}

A proof of the following result, probably folklore, can be found 
in~\cite{lytchak-geom-resol}.

\begin{prop}
  Let $G$ be a connected compact Lie group of isometries
of a simply-connected complete Riemannian manifold,  
consider the orbit space $X=M/G$, and the subset $X_{orb}$ 
of orbifold points of~$X$. Then $X_0:=X_{orb}\setminus\partial X_{orb}$ 
is exactly the set of non-singular (that is, principal and exceptional)
$G$-orbits. Furthermore $\pi_1^{orb}(X_0)=\{1\}$.
\end{prop} 

\Pf Due to Theorem~\ref{lt1}, $B=X_{orb}$ consists precisely
of the projections of points in $M$ where the slice
representation is polar. Let us show that all points
$x\in B\setminus\partial B$ represent non-singular $G$-orbits in $M$.
In fact, let $x\in B$ represent a singular orbit. Choose $p\in M$
projecting to $x$. Since the slice representation at~$p$ is polar,
for the normal slice $S$ at~$p$ we have that $S/G_p$ is isometric
to a Weyl chamber. The projection $S\to X$ induces an
open isometric embedding
$S/G_p\to B$. Since the Weyl chamber has non-empty boundary, so does
its image in $B$. Hence any neighborhood of~$x$ in $B$
contains boundary points. Since the boundary is closed, we deduce that
$x\in\partial B$. 

Under the assumption that $M$ is simply-connected, we now show that
all points in~$\partial B$ represent singular orbits. 
Indeed due to Remark~\ref{rem1}, the subset of non-singular orbits
$B_0\subset B$ has $\pi_1^{orb}(B_0)=\{1\}$, so~$B_0$ cannot contain strata
of codimension~$1$. But $B_0$ is open,
so, if it has a point in~$\partial B$, then it has a point
lying in a stratum of
codimension~$1$, and hence the whole stratum is contained in $B_0$, a
contradiction. We have shown that $X_0=B_0$, as desired. \EPf 

\medskip

We mention two further instances of major applications 
of Theorem~\ref{lt1}. 
A polar action is called \emph{infinitesimally polar} if 
all of its slice representations are polar. Due to 
Proposition~\ref{slice-representation-polar}, every polar 
representation is infinitesimally polar. 
Theorem~\ref{lt1} was the main tool in~\cite{gorodski-lytchak}
to classify the infinitesimally polar actions 
on Euclidean spheres. In particular:

\begin{thm}[Gorodski-Lytchak]
An isometric quotient $X$ of the unit sphere by a compact Lie group
is isometric to a Riemannian orbifold if and only if the 
universal orbi-covering $\tilde X$ of~$X$ 
is a weighted complex or quaternionic
projective space, or $\tilde X$ has constant curvature~$1$ or~$4$.
\end{thm}

Further, in~\cite{gorodski-kollross} this classification 
of quotients isometric to Riemannian orbifolds 
was extended to compact rank one symmetric spaces.

\section{Lecture 3: Variationally complete actions}

In the 1950s Bott and Samelson introduced the concept of variatonally
completess as a means of studying the topology of symmetric spaces and
their loop spaces~\cite{bott-samelson}.
Roughly speaking, a proper and isometric action 
on a complete Riemannian manifold is variationally complete if it
produces enough Jacobi fields along geodesics to determine the multiplicities
of the focal points to the orbits. In modern terminology,
they proved that the orbits of a variationally
complete action are taut submanifolds of the ambient space,
in the sense that the energy functional on the space of curves\ft{
  Of $H^1$-Sobolev class.} joining a generic point to an
arbitrary, fixed orbit is a perfect Morse function. This establishes
strong relations between the topology of the ambient manifold
and the topology of the orbits~(cf.~\cite{terng-thorbergsson2}). 

\subsection{Bott and Samelson's and related results}

Let a Lie group $G$ act properly and isometrically on a complete
Riemannian manifold~$M$. 
A geodesic $\gamma$ of $M$ is called \emph{horizontal} if it is orthogonal
to one orbit (and hence to every orbit it meets). Fix a $G$-orbit $N$
and a horizontal geodesic $\gamma$ meeting $N$ at time $t=0$.
A Jacobi field along $\gamma$ is called an \emph{$N$-Jacobi field}
if it is the variational field of a variation of $\gamma$
through horizontal geodesics starting at $N$. Finally, the action $(G,M)$ is
called \emph{variationally complete} if for every $G$-orbit~$N$ and
every horizontal geodesic~$\gamma$ starting at~$N$, every
$N$-Jacobi field that vanishes for some $t>0$ is the restriction
of a $G$-Killing vector field along~$\gamma$.\ft{It is equivalent
  to require that every $N$-Jacobi field that is tangent to
another orbit is the restriction of a $G$-Killing vector field.} 

The motivation of Bott and Samelson to consider variationally
complete actions of $G$ on $M$ 
was to construct an explicit
basis of cycles in the $\Z_2$-homology of the path space $\Omega(M;x,N)$,
where $N$ is an arbitrary $G$-orbit, $x\in M$, and the paths start at~$x$ and
end at a point in $N$. In modern terminology, we can state their
result as follows: 

\begin{thm}[Bott-Samelson]
  The orbits of a variationally complete action are taut submanifolds
  (with respect to $\Z_2$-coefficients).
\end{thm}

Here a submanifold $N$ of $M$ is called \emph{taut} if, for every
nonfocal point $x$, the energy functional $E:\Omega(M;x,N)\to\R$,
$E(\gamma)=\frac12\int||\gamma'||^2ds$, is a perfect Morse function,
that is, every critical point (geodesic) of $E$ corresponds to a basis
element of $H_*(\Omega(M;x,N))$. Indeed, Bott and Samelson provide
an algorithm to construct an explicit cycle for each critical point.
In the same paper, for a symmetric space $G/K$,
they prove that the isotropy action of $K$ on $G/K$, the $K\times K$-action
on $G$ by left and right-multiplication, and the linear isotropy
action of $K$ on $T_{x_0}(G/K)\cong\Lp$ are variationally complete.
Soon thereafter, Hermann~\cite{hermann} found a more general family of
variationally complete actions on symmetric
spaces. Namely, if $K$ and $H$ are both symmetric subgroups of the
compact Lie group $G$, then the action of $H$ on $G/K$ is
variationally complete.

L. Conlon was a student of Conlon. 
In~\cite{conlon} he proved the following theorem:

\begin{thm}[Conlon]
A hyperpolar action of a compact Lie group $G$
on a complete Riemannian manifold $M$ is 
variationally complete.
\end{thm}

\Pf Let $N=Gp$ be a fixed
orbit and let $q$ be a focal point of $N$ (that is, a critical
value of the normal exponential map) along a geodesic
$\gamma:[0,\ell]\to M$ with $\gamma(0)=p$ and $\gamma(\ell)=q$.
Then there exists a Jacobi field $J$ along $\gamma$
satisfying $J(0)\in T_pN$,
$J'(0)+A_{\gamma'(0)}J(0)\in\nu_pN$
and $J(\ell)=0$; denote by $V$ the space of Jacobi fields
satisfying the first two of these conditions, and note that 
$\dim V=\dim M$. 

Fix $s_0\in(0,\ell)$ such that $p_0=\gamma(s_0)$ is a regular point 
for the action of $G$ and $p_0$ is not a focal point of $N$. 
There exists a unique section $\Sigma$ passing through
$p_0$. Of course, $\Sigma$ is flat and contains the image of $\gamma$. 
Since $p_0$ is not a focal point of $N$, the map
$J\in V\mapsto J(s_0)\in T_{p_0}M$ is a linear isomorphism. 

Decompose $J=J^v+J^h$ where $J^h$ is the orthogonal projection of $J$
on $\Sigma$. Due to the total-geodesicness of $\Sigma$, both 
$J^v$ and $J^h$ are Jacobi fields along $\gamma$. Since 
$J^h$ vanishes at $s=0$ and $s=\ell$ and $\Sigma$ is flat, we have
$J^h\equiv0$. Since $p_0$ is a regular point, $J^v(s_0)\in T_{p_0}(Gp_0)$. 
Let $X\in\mathfrak g$ be such that $X^*_{p_0}=J^v(s_0)$. 
Owing to $X^*\circ\gamma\in V$, we have
$X^*\circ\gamma=J^v=J$, finishing the proof. \EPf  

\subsection{The converse results}

It was proved in~\cite{gorodski-thorbergsson-repr}, by means
of classification, that
a variationally complete representation is orbit-equivalent
to the isotropy representation of a symmetric space, and hence is polar.
In~\cite{discala-olmos}, a direct proof
of this result was provided. Since the idea of the proof is
very simple and geometric, we present it in the sequel. 

\begin{thm}[Di Scala-Olmos]\label{thm:discala-olmos}
  A variationally complete representation of a compact Lie
  group $G$ on an Euclidean space~$V$ is polar.
\end{thm}

\Pf Let $p\in V$ be a regular point so that $N=Gp$ is a
principal orbit. Owing to Lemma~\ref{exponential-normal-space},
$\Sigma:=\nu_pN$ meets
all orbits.

Choose $\xi\in \nu_pN$ such that the
Weingarten operator $A_\xi$ has all eigenvalues nonzero.
This is possible, since $A_p=-\mathrm{id}$, and indeed
the set of such vectors is open and dense in $\nu_pN$.
Consider the geodesic $\gamma(s)=p+s\xi$, normal to $N$,
and fix $s_1>0$ such that $N_1=Gq$, $q=\gamma(s_1)$, is also a principal
orbit. Due to the assumption of variational completeness,
$q$ is not a focal point of $N$ along $\gamma$. We will
show that $T_pN=T_qN_1$ as subspaces of $V$. 

Each eigenvector $u\in T_pN$ of $A_v$, with corresponding
eigenvalue $\lambda\neq0$, gives rise to a Jacobi field
$J(s)=(1-\lambda s)u$
along the geodesic $\gamma(s)=p+s\xi$, associated to the variation
$\gamma_t(s)=c(t)+s\hat\xi(t)$, where $c$ is a smooth curve
in $N$ with $c(0)=p$ and $c'(0)=u$, and $\hat\xi$ is the parallel
extension of $\xi$ to a normal vector field along $c$. Since
$J(0)=u\in T_pN$ and $J(\frac1\lambda)=0$, the assumption
of variational completeness yields a Killing vector field $X$ induced
by $G$ such that $X\circ \gamma= J$. In particular,
$J(s)\in T_{\gamma(s)}(G\gamma(s))$ for all $s$.
Since $q$ is not a focal point of $N$ along $\gamma$,
$s_1\neq\frac1\lambda$ and therefore $u\in T_qN_1$. 
As the eigenvectors of $A_\xi$ span $T_pN$, this shows
$T_qN_1=T_pN$.

We have seen that $\Sigma$ is orthogonal to all principal
orbits passing through an open and dense subset of itself.
By a continuity argument, $\Sigma$ is orthogonal to every orbit it meets.
This finishes the proof. \EPf

\medskip

An isometric action of a compact Lie group on a compact symmetric
space can be lifted to a proper and Fredholm action of a
Hilbert-Lie group on a Hilbert space via the so-called
``holonomy map'', see~\cite{terng-thorbergsson}.
This idea was combined with 
the basic idea of the proof of Theorem~\ref{thm:discala-olmos}
to prove the following result in~\cite{gorodski-thorbergsson-varcomp}: 

\begin{thm}[Gorodski-Thorbergsson]\label{thm:gorodski-thorbergsson}
  A variationally complete action of a compact Lie group on a
  compact symmetric space is hyperpolar.
\end{thm}

The following result was proved in~\cite{lytchak-thorbergsson-varcomp}
and generalizes Theorems~\ref{thm:discala-olmos}
and~\ref{thm:gorodski-thorbergsson}. Its proof is the main goal 
of this lecture.

\begin{thm}[Lytchak-Thorbergsson]\label{thm:lytchak-thorbergsson-varcomp}
A variationally complete action on a complete Riemannian
manifold with nonnegative sectional curvature is hyperpolar.
\end{thm}

A special case of this theorem is when the group is discrete. It says
that a complete Riemannian manifold without conjugate points and with
nonnegative sectional curvature is flat. 
We shall explain the proof of Theorem~\ref{thm:lytchak-thorbergsson-varcomp}
in the sequel. The main tool is Wilking's transversal Jacobi
equation~\cite{wilking}, which can be viewed as an extension of the
methods of Morse theory from the case of Riemannian submersions
to the case of singular Riemannian foliations and, in particular,
isometric actions.

Let $\pi:M\to B$ be a Riemannian submersion
with horizontal and vertical distributions $\mathcal H$ and
$\mathcal V$, respectively. Then one of O'Neill's
equation says that the sectional curvature of a 
horizontal $2$-plane $\sigma\subset\mathcal H$
and its projection $\sigma^*=d\pi(\sigma)\subset TB$  are related
by
\begin{equation}\label{oneill}
 K(\sigma^*)=K(\sigma)+ 3||A_XY||^2, 
\end{equation}
where $A:\mathcal H\times\mathcal H\to\mathcal V$ is one of
O'Neill's tensors associated to $\pi$, namely,
\[ A_XY = (\nabla_{X^h}Y^h)^v+(\nabla_{X^h}Y^v)^h, \]
for all $X$, $Y\in TM$. The following properties of $A$ are easily
established:
\begin{enumerate}
\item $A_X{\mathcal H}\subset\mathcal V$ and 
$A_X{\mathcal V}\subset\mathcal H$ for all $X\in TM$. 
\item $A_X$ is skew-symmetric on $T_pM$ for all $p\in M$. 
\item $A_XY=\frac12[X,Y]^v$ for all $X$, $Y\in\mathcal H$. 
\end{enumerate}
In particular $\mathcal H$ is integrable if and only if 
$A_X\mathcal H=0$ for all $X\in\mathcal H$. 

Let now $X$, $Y$ be an orthonormal basis of the horizontal $2$-plane~$\sigma$.  
Then the right hand-side of O'Neill's equation~(\ref{oneill}) reads
\[ \langle R(Y,X)X,Y\rangle + 3\langle A_XY,A_XY\rangle 
= \langle (R(Y,X)X)^h- 3A^2_XY,Y\rangle.  \]
Every projectable Jacobi field $J$ along a horizontal geodesic~$\gamma$ 
projects to a Jacobi field along $\bar\gamma=\pi\circ\gamma$, and 
this projection induces an isomorphism between space of  projectable 
Jacobi fields modulo the vertical Jacobi fields along $\gamma$,
and the space of Jacobi fields along $\bar\gamma$. It follows 
that for every projectable Jacobi field $J$ along $\gamma$, the 
horizontal component $J^h$ satifies the ``transversal Jacobi equation''
\[ (J^h)''+ (R(J^h,\gamma')\gamma')^h - 3A^2_{\gamma'}J^h=0. \]

Let now $(G,M)$ be a proper and isometric action. On the regular
part there is a Riemannian submersion $\pi:M_{reg}\to M_{reg}/G$, 
to which the above considerations apply. Moreover, Wilking explained 
how to overcome the difficulties that arise when the horizontal 
geodesic $\gamma$ passes through singular points of the action;
here we assume $\gamma$ is a regular complete geodesic in the sense
that it passes through regular points; it follows that the
singular points along~$\gamma$ form a discrete set of parameters. 

The first step is to extend the principal
horizontal distribution $\mathcal H$ along the regular part
of the horizontal geodesic $\gamma$ 
to a smooth distribution defined everywhere along $\gamma$, namely
$\mathcal H_t\subset T_{\gamma(t)}M$ for all $t\in\R$, where
$\mathcal H_t=(T_{\gamma(t)}(G\gamma(t)))^\perp$ if $\gamma(t)$ is a
regular point of the $G$-action.
Let $\Lambda$ be the space of $N_t$-Jacobi fields along~$\gamma$,
where $N_t=G\gamma(t)$, and let $\Upsilon$ be the subspace of
vertical Jacobi fields. We put $\mathcal H_t:=(\mathcal V_t)^\perp$,
where we define
\[ \mathcal V_t:=\{J(t)|J\in\Upsilon\}\oplus\{J'(t)|J\in\Upsilon,\
  J(t)=0\}. \]

\begin{prop}
  $\mathcal V_t$ is a smooth vector bundle 
  of rank $\dim\Upsilon$ along $\gamma$. 
Moreover, $\mathcal V_t=T_{\gamma(t)}N_t$ if
  $\gamma(t)$ is regular.
\end{prop}

\Pf It is easy to see that $\{J(t)|J\in\Upsilon\}=T_{\gamma(t)}N_t$
for all~$t$. Indeed denote the Lie algebra of~$G$ by~$\Lg$. 
Then $T_{\gamma(t)}N_t=\mathrm{span}\{X^*_{\gamma(t)}|X\in\Lg\}$
and, for each $X\in\Lg$, the induced Killing field $X^*$ on $M$ 
restricts to a Jacobi field $J$ along $\gamma$ which 
is tangent to the $G$-orbits everywhere, that is, $J\in\Upsilon$. 

Suppose now $\gamma(t_0)$ is regular for some~$t_0$. Then 
$\{J'(t_0)|J\in\Upsilon,\ J(t_0)=0\}$ is trivial. 
Indeed if there was $J\in\Upsilon$ 
with $J(t_0)=0$ and $J'(t_0)\neq0$, then $\gamma(t_0)$  would be a 
focal point of all $N_t$, but this is impossible as $N_{t_0}$ is a 
principal orbit. 

Finally, if $Y$ is a smooth vector field along $\gamma$ that 
has an isolated zero at $t_1$, then the vector field $\tilde Y$, given by
\[ \tilde Y(t):=\left\{\begin{array}{ll}
                    \frac1{t-t_1}Y(t), &\mbox{if $t\neq t_1$,}\\
                    Y'(t_1), &\mbox{if $t=t_1$,} \end{array}\right. \]
is smooth, and the span of $Y(t)$, $\tilde Y(t)$ is one-dimensional 
on a neighborhood of~$t_1$. This proves all statements. \EPf   

\medskip

The second step is to extend the definition of~$A_{\gamma'(t)}$
to all $t$, namely, a skew-symmetric 
operator $A_t$ on $T_{\gamma(t)}M$ such that 
$A_t$ coincides with the O'Neill tensor $A_{\gamma(t)}$ if~$\gamma(t)$ is 
regular. Let $Y$ be a smooth vector field along~$\gamma$. Set
\[ A_tY(t) := ((Y^h)'(t))^v+((Y^v)'(t))^h. \] 
The tensor $A_t$ clearly satisfies the requirements.  

The last step is to write the differential equation along~$\gamma$ 
that vector fields of the form $Y=J^h$ for some $J\in\Lambda$ 
must satisfy. This equation was derived in~\cite{wilking}, and reads
\begin{equation}\label{wte}
 \frac{(\nabla^h)^2}{dt^2}Y  + (R(Y,\gamma')\gamma')^h-3A_t^2Y =0. 
\end{equation}
where $\frac{\nabla^h}{dt}$ is the connection induced on horizontal vector 
fields along $\gamma$ (see subsection~\ref{sec:wtje}).

The idea of the proof of Theorem~\ref{thm:lytchak-thorbergsson-varcomp} 
is roughly as follows. Owing to Theorem~\ref{thm:hlo},
it suffices to show that $\mathcal H$ is integrable over $M_{reg}$.
This, in turn, follows if $A_t$ vanishes identically. On the other hand,
if $A_t$ is not identically zero, due to nonnegative curvature of $M$,
we obtain that $M/G$ ``has positive curvature'' somwhere and thus
``conjugate points''. But variationally completeness of $(G,M)$
more or less means ``absence of conjugate points'' in $M/G$, leading to a
contradiction. 

For each $t$ the operator
\[ \mathcal R(t):v\mapsto (R(v,\gamma'(t))\gamma'(t))^h-3A_t^2v \]
is self-adjoint and positive semidefinite on $\mathcal H_t$. 
Therefore
\begin{equation}\label{ode}
 Y''(t)+\mathcal R(t)Y(t) = 0 
\end{equation}
is a differential equation of Morse-Sturm type (generalization
of the Jacobi equation), where the ``prime'' refers to $\frac{\nabla^h}{dt}$,
to which is associated an index form
\[ I_{a,b}(X,Y)=\int_a^b\langle X',Y'\rangle - \langle \mathcal R(t)X,Y\rangle\, dt,
\]
for each $a<b$, 
where $X$ and $Y$ are piecewise smooth sections of the horizontal distribution
along $\gamma|_{[a,b]}$ vanishing at~$a$ and $b$.
Suppose, to the contrary, that there is $t_0$ such that
$A_{t_0}\neq0$. Then there is a unit vector $v_0\in\mathcal H_{t_0}$ such that
$\langle \mathcal R(t_0)v_0,v_0\rangle>0$. Let $Z_0$ be the 
$\nabla^h$-parallel vector field
 along $\gamma$ such that $Z_0(t_0)=v_0$. Then
\begin{equation}\label{C}
 C := \int_{t_0-1}^{t_0+1}\langle \mathcal R(t)Z_0(t),Z_0(t)\rangle\, dt>0. 
\end{equation}
Let $\varphi:\R\to[0,1]$ be a (smooth) bump function with support
contained in $[t_0-N,t_0+N]$ and satisfying $\varphi|_{[t_0-1,t_0+1]}\equiv1$,
for some $N>1$. In fact, we can take $N$ as large as needed to further
have   
\[ \int_{t_0-N}^{t_0+N}\varphi'(t)^2\,dt<C. \]
Set $Z(t)=\varphi(t)Z_0(t)$. For $a=t_0-N$ and $b=t_0+N$,
our choices yield
\begin{align*}
  I_{a,b}(Z,Z) & = \int_{t_0-N}^{t_0+N}||Z'(t)||^2-\langle \mathcal R(t)Z(t),Z(t)\rangle\, dt \\
  &=\int_{t_0-N}^{t_0+N}\varphi'(t)^2\,dt-\int_{t_0-N}^{t_0+N}\varphi(t)^2\langle \mathcal R(t)Z_0(t),Z_0(t)\rangle\, dt \\
  &< C - \int_{t_0-1}^{t_0+1}\langle \mathcal R(t)Z_0(t),Z_0(t)\rangle\, dt \\
&=0.
\end{align*}
This shows that $I_{a,b}$ has negative index. By the Sturm 
oscillation theorem, which is a generalization of the Morse index theorem, 
there is a nonzero solution~$Y$ of~(\ref{ode}) such that 
$Y(a)=Y(c)=0$ for some $c\in(a,b)$ (a ``conjugate point'').  
By continuity of the index, if necessary, we may perturb $a$ slightly so 
that $\gamma(a)$ and $\gamma(c)$ become regular points of the $G$-action
and remain conjugate points. 

Finally, we show that $Y$ has the form $Y=J^h$ for some $J\in\Lambda$. 
There is an $N_a$-Jacobi field~$J$ along $\gamma$ with initial conditions
$J(a)=0$ and $J'(a)=\frac{\nabla^hY}{dt}(a)$. Set $\hat Y=J^h$. 
Then $Y$ and $\hat Y$ are both solutions of the same differential 
equation~(\ref{wte}), with the same initial conditions at~$a$. 
It follows that they coincide. This yields a Jacobi field~$J$
along~$\gamma$, tangent to orbits at $t=a$ and $t=c$, which is 
not tangent to orbits everywhere, leading to a contradiction 
with the assuption that $(G,M)$ is a variationally complete action.

We have shown that $A_t\equiv0$ for all~$t$, which says that 
$\mathcal H$ is integrable over~$M_{reg}$, and this implies that 
its integral manifolds can be extended to sections of $(G,M)$,
due to Theorem~\ref{thm:hlo}. It only remains to check 
that sections are flat. Let $\Sigma$ be a section of $(G.M)$.
Since $\Sigma$ is totally geodesic, it is nonnegatively curved. 
Suppose, to the contrary, that $\Sigma$ has a tangent~$2$-plane~$\sigma$ 
with positive 
curvature at a point~$p\in\Sigma$, which we can assume is a 
regular point of the $G$-action, by denseness. Let $\gamma$ 
be a horizontal geodesic with $\gamma(0)=p$ and~$\gamma'(0)\in\sigma$. 
Then we can find $v_0\in\mathcal H_0$ such that $\langle\mathcal R(0)v_0,v_0\rangle>0$,
and proceed as in~(\ref{C}) to reach a contradiction. Hence~$\Sigma$ 
is flat and this completes the proof of
Theorem~\ref{thm:lytchak-thorbergsson-varcomp}. 

\section{Appendix}

\subsection{An algebraic
  criterion for polar actions on symmetric spaces}\label{criterion}

The following criterion was proved in~\cite{gorodski-polar}.

\begin{prop}[Gorodski]
Let $M=G/K$ be a symmetric space without Euclidean factor  endowed 
with a Riemannian metric induced from some $\mbox{Ad(G)}$-invariant
inner product on the Lie algebra $\Lg$ of $G$. 
Consider a closed, connected subgroup $H\subset G$.
By replacing $H$ by a 
conjugate, if necessary, we may assume that $\bar1=1K\in G/K$ is a regular 
point. Write $\Lg=\Lk+\Lp$ for the Cartan decomposition, denote by $\Lh$ 
the Lie algebra of $H$, and define 
$\Lm=\Lp\cap\Lh^\perp$. Then the action of $H$ on $M$ is polar if and only if
the following two conditions hold:
\begin{enumerate}
\item[(i)] $\Lm$ is a Lie triple system (LTS), 
that is $[\Lm,[\Lm,\Lm]]\subset
\Lm$; and
\item[(ii)] $[\Lm,\Lm]\perp\Lh$. 
\end{enumerate}
\end{prop}

\Pf Let $\pi:G\to G/K$ be the canonical projection, and
for $a\in G$ write $\bar{a}=\pi(a)=aK\in G/K$. We have that
$H\times K$ acts on $G$ by left and right translations, 
$H$ acts on $G/K$ by left translations,
and $\pi$ is an equivariant Riemannian submersion. 
Therefore the tangent space to the orbit $H\cdot\bar a$ 
at the point $\bar a$ is
\[ T_{\bar a}(H\cdot\bar a)=d\pi_a(\Lh\cdot a+a\cdot\Lk), \]
so that 
\[ a^{-1}\cdot T_{\bar a}(H\cdot\bar a)=\pi_*(\Ad_{a^{-1}}\Lh) \]
where $\pi_*:\Lg\to\Lp$ is the projection. Taking orthogonal
complements we get that 
\[ a^{-1}\cdot \nu_{\bar a}(H\cdot\bar a)=\Lp\cap\Ad_{a^{-1}}\Lh^\perp \]
where $\nu_{\bar a}(H\cdot\bar a)$ denotes the normal 
space to $H\cdot\bar a$ 
at $\bar a$. In particular we have
\[ \nu_{\bar 1}(H\cdot\bar1)=\Lm. \]
Since $\bar 1$ is a regular point, the action of $H$ on $M$ is polar
if and only if $\Sigma=\mathrm{Exp}_{\bar1}[\Lm]$ is a section,
where $\mathrm{Exp}$ denotes the Riemannian exponential map of $M$. 

If $\Sigma$ is a section, then it is totally geodesic in $M$.
Now the fact that $\Lm=T_{\bar1}\Sigma$ is a LTS is a standard fact from
the theory of symmetric spaces, see~\cite{helgason}, p.~224. 
This proves condition (i)
and shows that $\Ls=[\Lm,\Lm]+\Lm$ is a subalgebra
of $\Lg$. Let $S$ denote the corresponding connected subgroup of $G$.
Now $\pi(S)=\Sigma$ and the elements of $S$ induce isometries of
$\Sigma$. Let $a\in S$. Since $\Sigma$ intersects the orbits of $H$
orthogonally, we have that 
$T_{\bar a}\Sigma\subset \nu_{\bar a}(H\cdot\bar a)$ and therefore
\[ \Lm=T_{\bar1}\Sigma=a^{-1}\cdot T_{\bar a}\Sigma
\subset a^{-1}\cdot \nu_{\bar a}(H\cdot\bar a) = 
\Lp\cap\Ad_{a^{-1}}\Lh^\perp. \]
This proves that $\Ad_{a}\Lm\subset\Lh^\perp$, so by taking 
$X$, $Y\in\Lm$ arbitrary and $a=\exp tX$ we get that
$\Ad_{\exp tX}Y\in\Lh^\perp$ and hence
$[X,Y]=\frac{d}{dt}\big|_{t=0}\Ad_{\exp tX}Y\in\Lh^\perp$.
This gives condition (ii) and proves half the proposition.

Conversely, if conditions (i) and (ii) hold, and $\Ls$, $S$ 
are as above, then $\Ls\subset\Lh^\perp$
so $\Ls=\Ad_{a^{-1}}\Ls\subset\Ad_{a^{-1}}\Lh^\perp$ for $a\in S$
and then $\Lm=\Ls\cap\Lp\subset\Ad_{a^{-1}}\Lh^\perp\cap\Lp$.
This implies that 
$T_{\bar a}\Sigma\subset \nu_{\bar a}(H\cdot\bar a)$
for $\bar a\in\Sigma$ and hence $\Sigma$ intersects 
the orbits of $H$ orthogonally. The fact that $\Sigma$ intersects
\emph{all} the orbits of $H$ follows from 
Lemma~\ref{exponential-normal-space}, and this 
finishes the proof of the proposition. \EPf

\begin{cor}\label{hyperpolar}
The action of $H$ on $M$ is hyperpolar if and only if $\Lm$ is an
Abelian subalgebra of~$\Lp$.
\end{cor}

\Pf If $\Lm$ is an Abelian subalgebra of $\Lp$ then the criterion 
of the proposition is satisfied, and $\Sigma$ is flat by the curvature  
formula for symmetric spaces. Conversely, if the action of $H$ on $M$ 
is hyperpolar, then $\Sigma$ is flat so $[[\Lm,\Lm],\Lm]=0$. 
But
\[ \langle [\Lm,\Lm],[\Lm,\Lm]\rangle = \langle [[\Lm,\Lm],\Lm],\Lm\rangle
=0. \]
Since $[\Lm,\Lm]\subset\Lk$
and the Killing form of $\Lg$ is negative definite on $\Lk$,
we deduce that $\Lm$ is Abelian, as wished. \EPf

\begin{cor}[Hermann actions]
  If $H$ is a symmetric subgroup of $G$, then its action on
  $M$ is hyperpolar.
  \end{cor}

\Pf By replacing $H$ by a conjugate, we may assume that 
$\bar 1$ is a regular point. 
Write $\Lg=\Lh+\Lq$ for the involutive decomposition 
realtive to the symmetric pair $(G,H)$. Then $\Lm$ as in the 
proposition is given by $\Lp\cap\Lq$, and we only need to see that it is 
Abelian. 

The orbit through $\bar 1$ is $HK/K=H/K\cap H$, and 
$\nu_{\bar 1}(HK/K)=\Lm$. Since $\bar1$ is a regular point,
the slice representation $(K\cap H,\Lm)$ is trivial. 
In particular $[\Lk\cap\Lh,\Lm]=0$. Since $[\Lm,\Lm]\subset\Lk\cap\Lh$,
we see that $[[\Lm,\Lm],\Lm]=0$, and hence $[\Lm,\Lm]=0$ as 
in Corollary~\ref{hyperpolar}. \EPf

\subsection{Cartan's and Hermann's theorems}\label{cartan-hermann}

Let $M$ be a Riemannian manifold. 
In his book on Riemannian geometry,
Cartan states a criterion for the local existence
of a totally geodesic submanifold in~$M$ with a given tangent space
at a given point.

\begin{thm}[Cartan]\label{thm:cartan}
   Let $M$ be a Riemannian manifold, 
fix a point $p\in M$ and a subspace $S\subset T_pM$.
Assume that there is a normal ball \[ V=\exp_p(B(0_p,\epsilon)) \]
such that for every unit speed
radial geodesic $\gamma:[0,\ell]\to V$
emanating from $p$ ($\ell<\epsilon$), 
\begin{equation}\label{condition2}
 R(P_\gamma(u),P_\gamma(v))P_\gamma(w)\in P_\gamma(S), 
\end{equation}
 for every  $u$, $v$, $w\in S$, where $P_\gamma$ denotes the parallel 
transport along~$\gamma$, from $0$ to~$\ell$.  Then there exists a 
totally geodesic submanifold manifold~$N$ of~$M$ such that $T_pN=S$. 
\end{thm}

\Pf Let $N=\exp_p(S\cap B(0_p,\epsilon))$. We will explain
why $N$ is totally geodesic. It suffices to see that parallel
transport in $M$, along piecewise smooth curves in $N$, preserves the tangent
spaces of $N$. 

In the case of a radial geodesic $\gamma(t)=\exp_p(tv)$ with~$v\in S$ and
$||v||<\epsilon$,
$t\in[0,1]$, this follows from the Jacobi equation.
In fact consider $q=\gamma(t_0)$ for some $t_0\in(0,1)$.
Recall that the Jacobi field $J$ along $\gamma$ with $J(0)=0$ and
$J'(0)=u\in S$ is given by $J(t)=d(\exp_p)_{t_0v}(t_0u)$. 
Let $E_1=\gamma'$, $E_2$,\ldots,$E_n$ be a parallel
orthonormal frame along~$\gamma$, where $E_1(0)\ldots, E_k(0)$ are tangent
to $S$ and  $E_{k+1}(0)\ldots, E_n(0)$ are normal to $S$. 
Write $J=\sum_i a_iE_i$. Then $-a_i''+\sum_j\langle R(E_1,E_j)E_1,E_i\rangle a_j = 0$ for all~$i$, where $a_i(0)=0$ for all $i$, and $a_i'(0)=0$ for $i>k$.
Owing to~(\ref{condition2}),
$\langle R(E_1,E_j)E_1,E_i\rangle\equiv0$ for $i>k$ and $j\leq k$,
so we deduce that $a_i\equiv0$  vanishes identically for $i>k$. Therefore
$J$ is everywhere tangent to the parallel translates of $S$ along $\gamma$.
Since $T_qN=d(\exp_p)_{t_0v}(S)$, this proves that the tangent spaces
of $N$ along $\gamma$ are parallel along $\gamma$.

In the case of an arbitrary piecewise smooth curve $\eta:[0,1]\to V$,
we join each point $\eta(s)$ to~$p$ by a radial geodesic, so as
to obtain a parametrized surface $f(s,t)$, $(s,t)\in[0,1]\times[0,1]$,
where $\gamma_s=f(s,\cdot)$ is a radial geodesic for each~$s$, and
$f(s,0)=p$, $f(s,1)=\eta(s)$ for all~$s$. Consider the vector fields
along $f$ given by
\[ \frac{\bar\partial}{\partial s}=f_*\frac{\partial}{\partial s}
  \quad\mbox{and}\quad
  \frac{\bar\partial}{\partial t}=f_*\frac{\partial}{\partial t}. \]
Then  $\frac{\bar\partial}{\partial s}$ is a Jacobi field along
each~$\gamma_s$, and we already know from the argument in the
previous paragraph that
it is everywhere tangent to~$N$ and sits in the parallel translate
of~$S$ along $\gamma_s$, whose value at~$f(s,t)$ we denote by $S_{s,t}$;
note that $T_{f(s,t)}N=S_{s,t}$. Further, it is clear that
$\frac{\bar\partial}{\partial t}$ is everywhere tangent to~$N$.
Let $z\in T_qN=S_{0,1}$ be arbitrary, where 
$q=\eta(0)=f(0,1)$.
We parallel translate~$z$ along $\gamma_0$ from $q$ to $p$ to obtain
$z_0\in S$,
and then we parallel translate~$z_0$ along each $\gamma_s$ from $p$
to $\eta(s)$ to obtain a vector field $Z$ along $f$. By the previous 
paragraph we know that $Z(s,t)\in S_{s,t}$. The main
calculation now is
\begin{equation}\label{R}
  R\left(\frac{\bar\partial}{\partial t},\frac{\bar\partial}{\partial s}\right)Z
  =\frac{\nabla}{dt}\frac{\nabla}{ds}Z,
\end{equation}
since $Z$ is parallel along each~$\gamma_s$. Let $w\in T_pM$ be normal to~$S$
and extend it to a vector field~$W$ 
along $f$ and parallel along each $\gamma_s$. Then $W(s,t)\perp S_{s,t}$
for all~$(s,t)$.  
We take the inner product of~(\ref{R}) throughout with $W$. Since the 
left hand-side of~(\ref{R}) lies in~$S_{s,t}$ by~(\ref{condition2}), we obtain 
\[ 0 = \langle  \frac{\nabla}{dt}\frac{\nabla}{ds}Z,W\rangle 
= \frac{d}{dt}\langle \frac{\nabla}{ds}Z,W\rangle. \]
Now $\langle\frac{\nabla}{ds}Z,W\rangle$ is constant in~$t$. 
Since it vanishes at $t=0$, it must also vanish at $t=1$. This proves
that $\{S_{s,1}\}_{s\in[0,1]}$ is a parallel family of subspaces
along~$\eta$, that is the tangent spaces to~$N$ are parallel along $\eta$. \EPf

\medskip
   
We next introduce some terminology. Let $M$ be a Riemannian manifold, 
fix a point $p\in M$ and a subspace $S\subset T_pM$. A 
once-broken geodesic $\gamma:[0,\ell]\to M$, emanating from $p$ 
and broken at $t_0\in(0,\ell)$, is called \emph{$S$-admissible} 
if $\gamma'(0)\in S$, $\gamma'(t_0^+)$ lies in the parallel transport
of $S$ along $\gamma$ from $p$ to $\gamma(t_0)$, and $\gamma|_{[t_0,\ell]}$ 
sits in a convex neighborhood of $\gamma(t_0)$. 

\begin{thm}[Hermann]\label{hermann}
 Let $M$ be a Riemannian manifold, 
fix a point $p\in M$ and a subspace $S\subset T_pM$.
Assume that for every $S$-admissible
once-broken geodesic $\gamma:[0,\ell]\to M$
emanating from $p$, 
\begin{equation}\label{condition}
 R(P_\gamma(u),P_\gamma(v))P_\gamma(w)\in P_\gamma(S), 
\end{equation}
 for every  $u$, $v$, $w\in S$, where $P_\gamma$ denotes the parallel 
transport along~$\gamma$, from $0$ to~$\ell$.  Then there exists a 
complete totally geodesic isometric 
immersion of Riemannian manifold~$N$ into~$M$ with $p\in N$ and
$T_pN=S$.
\end{thm}

\Pf By Cartan's local existence theorem~\ref{thm:cartan}, there exists
a totally geodesic immersion of a Riemannian manifold $N$ into a normal
neighboorhood of~$p$ such that $p\in N$ and $T_pN=S$. We take $N$
maximal and assume, by contradiction, that $N$ is not complete.
Then there is a geodesic $\gamma:[0,1)\to N$ such that
$\lim_{t\to 1-}\gamma(t)$ does not exist in $N$. By completeness of $M$,
$\gamma$ can be continued past $t=1$ to a complete geodesic
$\tilde\gamma:[0,+\infty)\to M$. Let $\tilde S$ be the parallel translate
of $S$ along $\tilde\gamma$ to $q=\tilde\gamma(1)$. 
Due to~\ref{condition} and using again Theorem~\ref{thm:cartan},
there is totally geodesic
submanifold $\tilde N$ with $q\in\tilde N$ and
$T_{\tilde q}\tilde N=\tilde S$. Of course the parallel translate of $S$
along $\tilde\gamma$ from $p$ to $\gamma(t_1)$ for $t_1\in(0,1)$
coincides with the parallel transport of $\tilde S$ along $\tilde\gamma$
from~$q$ to $\gamma(t_1)$. This shows that the tangent spaces of
$N$ and $\tilde N$ coincide at $\gamma(t)$
for $t=1-\delta$ with $\delta>0$ sufficiently
small, and hence implies that $N$ and $\tilde N$ coincide on a neighborhood
of $\gamma(0,1)$. Since $q\in\tilde N$, this contradicts the maximality
of~$N$. \EPf 

\subsection{Wilking's transversal Jacobi equation}\label{sec:wtje}

Let $(G,M)$ be a proper and isometric action, and 
let $\pi:M\to M/G$ be the projection.
In this appendix we reproduce the derivation of Wilking's transversal 
Jacobi equation~\cite{wilking} in this special case:
\begin{equation}\label{wte2}
 \left(\frac{\nabla^h}{dt}\right)^2Y + (R(Y,\gamma')\gamma')^h-3A_t^2Y =0, 
\end{equation}
along a horizontal geodesic $\gamma$ for vector fields of the form 
$Y=J^h$, where $J$ belongs to the space $\Lambda$ of 
$N_t$-Jacobi fields along $\gamma$, $N_t=G\gamma(t)$. 

We denote by $\mathcal{J}$ 
the space of all Jacobi fields 
along~$\gamma$. Consider the skew-symmetric bilinear
form $\omega(J_1,J_2)=\langle J_1',J_2\rangle
- \langle J_1,J_2'\rangle$, for $J_1$, $J_2\in\mathcal {J}$.
Then $\omega$ is constant and defines a symplectic form 
on $\mathcal{J}$ such that $\Lambda$ is a Lagrangian subspace
and $\Upsilon$ is an isotropic subspace. Indeed
an $N_t$-Jacobi field $J$ satifies that
$J(t)$ is vertical and $J'(t)+S_{\xi}J(t)$ 
is horizontal for this~$t$, 
where $\xi=\gamma'(t)$ and $S_\xi$ denotes the shape operator in the direction
of a normal vector $\xi$ to the orbit $N_t$, so   
\begin{align*}
\omega(J_1,J_2)&=\omega(J_1(t),J_2(t))\\
& =\langle J_1'(t),J_2(t)^v\rangle -\langle J_1(t)^v,J_2'(t)\rangle\\
         &=-\langle S_\xi J_1(t),J_2(t)\rangle +\langle J_1(t),S_\xi J_2(t)\rangle\\
&=0
\end{align*}
for $J_1$, $J_2\in\Lambda$; in addition, $\dim\Lambda=\frac12\dim\mathcal{J}$,
and $\Upsilon\subset\Lambda$. 

By continuity it suffices the check equation~(\ref{wte2}) at a point $t_0$
such that $N_{t_0}$ is a principal orbit. Since we can add a vertical
Jacobi field to~$J$ without changing~$Y$, we may assume that 
$J(t_0)\in\mathcal H_{t_0}$. Let $E_1,\ldots,E_n$ be a
$\nabla^h$-parallel orthonormal frame field of $\mathcal H$ 
along $\gamma$ with $E_1(t_0)=J(t_0)$. We first claim that 
\begin{equation}\label{claim1}
 (J'(t_0))^v = -A_{t_0} J(t_0).
\end{equation}
Indeed, for any vertical Jacobi field $V$ along $\gamma$, 
\begin{align*}
 \langle J'(t_0),V(t_0)\rangle &= \langle J(t_0),V'(t_0)\rangle\\
& = \langle J(t_0),(V'(t_0))^h\rangle\\
& = \langle J(t_0),A_{t_0}V(t_0)\rangle\\
&=-\langle A_{t_0}J(t_0),V(t_0)\rangle,
\end{align*}
where we have used that $\Lambda$ is Lagrangian and $A_{t_0}$ is 
skew-symmetric.

The second claim is that 
\begin{equation}\label{claim2}
E_i'(t)=AE_i(t) 
\end{equation}
for all~$t$. In fact, for any vertical Jacobi field $V$ along $\gamma$, 
\begin{align*}
\langle E_i'(t),V(t)\rangle & = -\langle E_i(t),V'(t)\rangle \\
& = - \langle E_i(t),(V'(t))^h\rangle \\
& = -\langle E_i(t),A_tV(t)\rangle \\
& = \langle A_tE_i(t),V(t)\rangle, 
\end{align*}
proving the claim.

Now we can finish the proof of~(\ref{wte2}) as follows.
First note that 
\begin{align*}
\langle E_1,E_k''\rangle_{t_0} & = -\langle E_1',E_k'\rangle_{t_0}+\frac{d}{dt}\Big|_{t=t_0}\langle E_1,E_k'\rangle\\
&=-\langle E_1',E_k'\rangle_{t_0},
\end{align*}
since $E_1$ is horizontal and $E_k'$ is vertical. Now
\begin{align*}
\langle\left(\frac{\nabla^h}{dt}\right)^2J^h,E_k\rangle & = 
\frac{d^2}{dt^2}\Big|_{t=t_0}\langle J^h,E_k\rangle \\
&=\frac{d^2}{dt^2}\Big|_{t=t_0}\langle J,E_k\rangle \\
&=\langle J'',E_k\rangle_{t_0}+2\langle J',E_k'\rangle_{t_0}+
\langle J,E_k''\rangle_{t_0} \\
&=-\langle R(J,\gamma')\gamma',E_k\rangle_{t_0}
+2\langle J',E_k'\rangle_{t_0}+
\langle E_1,E_k''\rangle_{t_0} \\
&=-\langle R(J^h,\gamma')\gamma',E_k\rangle_{t_0}
+2\langle J',E_k'\rangle_{t_0}-\langle E_1',E_k'\rangle_{t_0}\\
&=-\langle R(J^h,\gamma')\gamma',E_k\rangle_{t_0}
-2\langle A_{t_0}E_1,A_{t_0}E_k\rangle-\langle A_{t_0}E_1,A_{t_0}E_k\rangle\\
&=-\langle R(J^h,\gamma')\gamma',E_k\rangle_{t_0}
+3\langle A_{t_0}^2J^h,E_k\rangle,
\end{align*}
as wished.

\providecommand{\bysame}{\leavevmode\hbox to3em{\hrulefill}\thinspace}
\providecommand{\MR}{\relax\ifhmode\unskip\space\fi MR }
\providecommand{\MRhref}[2]{%
  \href{http://www.ams.org/mathscinet-getitem?mr=#1}{#2}
}
\providecommand{\href}[2]{#2}

\end{document}